\documentclass[leqno,twoside,11pt]{article}
\usepackage{amsfonts,amsmath}
\begin{document}
\newtheorem{lem}{Lemma}[section]
\newtheorem{defi}[lem]{Definition}
\newtheorem{theo}[lem]{Theorem}
\newtheorem{pro}[lem]{Proposition}
\newtheorem{remark}[lem]{Remark}
\newtheorem{cor}[lem]{Corollary}
\newcommand{\beqa}{\begin{eqnarray}}
\newcommand{\eeqa}{\end{eqnarray}}
\newcommand{\confe}{\stackrel{*}{\rightharpoonup}}
\newcommand{\R}{I\! \! R}
\newcommand{\N}{I\! \! N}
\newcommand{\Z}{\mathbb{Z}}
\newcommand{\proof}{{\bf Proof.-}}
\newcommand{\BO}{\sqcap \! \! \! \! \sqcup} 
\newcommand{\Rn}{I\! \! R^{n}}
\newcommand{\Le}{I\! \! L}
\newcommand{\fo}{\widehat}
\newcommand{\maxi}{| \! \! [}
\newcommand{\maxd}{] \! \! |}
\newcommand{\nor}{\bigg\|}
\newcommand{\NOR}{\big\bigg\|}
\setlength{\baselineskip}{1.3\baselineskip}
\title{\bf Asymptotic expansion for the models of nonlinear 
dispersive, dissipative equations}
\author{
Ra\'ul Prado\\
{{\footnotesize Departamento de Matem\'atica Universidade Federal de Paran\'a,}}\\ 
{{\footnotesize Caixa Postal 19081, CEP 81531-990 Curitiba, PR, Brazil}}\\
{{\footnotesize rprado@mat.ufpr.br} }\\
} 
\maketitle
\begin{abstract} 
Considered herein are the family of nonlinear equations with both
dispersive and dissipative homogeneous terms appended. Solutions of these 
equations that start with finite energia decay to zero as time goes to infinity.
We present an asymptotic form which renders explicit the influence of the
dissipative, dispersive and nonlinear effect in this decay.
We obtain the second term in the asymptotic expansion, as time goes to infinity,
of the solutions of this equations and the complete asymptotic expansion,  
as time goes to infinity, of the linearized equations.
\end{abstract}
\vspace*{5mm}
\renewcommand{\theequation}{\thesection.\arabic{equation}}
\setcounter{equation}{0}
%
\section{Introduction}
This paper is intended to study the asymptotic expansion of solutions of one 
family of nonlinear, dispersive equations under the effect of dissipation.
Our model equation takes the form 
\begin{equation}\label{gene}
\left\{ 
\begin{array}{ll}
u_t + M u_t + M u + u_x + \big( u^q\big)_x= 0, \quad x \in \R, \quad t>0;\\
\quad \; \;u(x,0) = u_0 (x),
\end{array}
\right.
\end{equation}
where $M$ is defined as Fourier multiplier homogeneus operator by 
\begin{gather}\label{opM}
\fo{Mf} (\xi) = |\xi|^{m} \fo{f}(\xi),\quad f \in H^{m}(\R) \quad m \geq 1,
\end{gather}
The circumflexes connote Fourier transform, and subscripts denote 
partial differentiation.
When equations of the class (\ref{gene}) arise as models of phisical phenomena
$u=u(x,t)$ represents the displacement of the medium of propagation 
from its equilibrium position and is a real-valued function of 
two real variables: $x$ (called the spatial variable) is proportional to 
distance in the direction of propagation 
and $t>0$ is proportional to time.

In this paper, $u^q$ schould be interpreted either as 
$|u|^q$ or $|u|^{q-1}\;u$.
We shall assume that $q>m>2$. All physical constant which may appear 
in (\ref{gene}) are put to be equal to $1$, for simplicity.

Equations of the form (\ref{gene}) 
arise when dissipation, dispersion, and the effect of nonlinearity 
are appended to the {\it transport equation} $u_t + u_x$ for the 
unidirectional wave propagation. 
The damping is represented here by $Mu$.
When $Mu=0$, this models arise in a wide variety of circunstances 
(see Biler \cite{Biler}, Benjamin \cite{Benja1,Benja2, Benja3} , 
Bona \cite{Bona1, Bona2} Abdelouhab et al \cite{ABFS} ).
The particular class in which the nonlinearity is a monomial,
the dispersive and dissipative terms are homogeneus, provides perhaps 
the simplest class of model in which to study the three effects.
The equation (\ref{gene}) is a more simple model of a general case 
than were studied by V. Bisognin and G. Perla in \cite{BP}, decay rates of 
the solutions in $\Le^{p}(\R)$ spaces, $2 \leq p \leq \infty$,
were obtained.

In the case $m=2$, i.e. when $Mu = - u_{xx}$, then the equation (\ref{gene})
is the well-known generalized Benjamin Bona Mahony Burger equation:
$$
u_t - u_{xxt} - u_{xx} + u_x + (u^q)_x= 0, \quad x \in \R, \quad t>0,
$$
this model appear when one attempt to describe the propagation of 
small-amplitude long waves in nonlinear dispersive media taking into 
account dissipative mechanisms. 
For solutions of this equations G. Karch obtain in \cite{Karch4} 
the first an second terms of the asymptotic expansion, when 
$t \rightarrow \infty$, of both linearized 
equation and  nonlinear equation. 

In \cite{razu}, we obtain the complete asymptotic expansion of solutions the 
linearized equation (the $n$-dimensional case), and compute the second term 
in the asymptotic expansion in the two-dimensional case, with quadratic 
nonlinear term.

When $m=4$ the equation (\ref{gene}) takes the form
\begin{gather*}
u_t + u_{xxxxt} + u_{xxxx} + u_x + (u^q)_x= 0, \quad x \in \R, \quad t>0,
\end{gather*}
this is the  Roseneau equation, with the term 
$u_{xxxx}$ associated with dissipative phenomena. The Cauchy
problem for this equation were solved by M. Park in \cite{Park}. The asymptotic 
expansion of solutions, when $t \rightarrow \infty$, is a consequence of our
results.

For the well-posedness of the initial-value problem (\ref{gene}), we refer to 
Bisognin and Perla \cite{BP}.
It is sufficient to know for the purpose of this paper that it is always 
possible to construct global in time solutions for any initial data 
$u_0 \in W^{2,1}(\R)$ provided either $\|u_0\|_{W^{2,1}(\R)}$ 
is small or some restriction on $q$ are imposed.
In the Section 6, we put this problem more carefully.

This paper is organized as follows. 
In the next Section we state and discuss main results
concerning to the equation (\ref{gene}). 
In Section 3 we prove preliminary result of the complete asymptotic 
expansion of solutions to the equation of type (\ref{gene}) linearized, 
we will use essentially the Taylor Theorem, the Plancherel equality, and 
the complete asymptotic expansion for heat equation of generalized type.
In section 4, we prove the complete asymptotic expansion of the linearized
equation (\ref{gene}) using essentially,
some properties of the Bessel Potential of order $\beta$, $K_{\beta}$.
Section 5 contains a result on the complete asymptotic expansion of 
generalized type KdV-B linear equation where the same techniques apply.
For completeness of the exposition, the global-in-time solutions to 
(\ref{gene}) are constructed  in Section 6. 
In Section 7 we calculate the second term in 
the asymptotic expansion, when $t \rightarrow \infty $, of the solution to 
the nonlinear equation (\ref{gene}). The proof bases on the same ideas as 
those in G. Karch \cite{Karch4}.\\
{\bf Notation.} The notation to be used is standard. 
For $1\leq p \leq \infty$, the $\Le^p (\R)$-norm of a Lebesgue measurable 
real-valued functions defined on $\Rn$ is denoted by $\| f \|_p$.
The Fourier transform of $u$ is given by  
$ {\cal{F}}u(\xi) = \fo{u}(\xi) \equiv \int_{\R} e^{-ix\cdot \xi}u(x) dx$. 
If $m \in \R$ we denote by $H^m (\R)$ the Sobolev space of order $m$ as the 
completion of the Schwartz space ${\cal{S}}(\R)$ respect to the norm
$\|u\|_{H^m} \equiv \biggl( \int_{\R}(1+|\xi|^{2})^m |\fo{u}(\xi)|^2 d\xi 
\biggr)^{1/2}$. 
For simplicity, we write $\int = \int_{\R}$.
The letter $C$ will denote generic positive constants, which do not depend on 
$u$, $x$ and $t$, but may vary from line to line during computations.
\section{Main result}
%
In order to eliminate the convective term
of order one in the equation (\ref{gene}), we define
the traslated function $v(x,t)=u(x-t,t)$. Then $u$ solves
(\ref{gene}) if and only if $v$ solves 
\begin{eqnarray}\label{gene1}
\left\{ \begin{array}{rcl}
v_{t} + Mv_{t} +Mv + Mv_x + ( v^{q} )_x & = & 0 
\:,\quad \quad \quad \quad  in \quad \R \times (0,\infty)\:; \\
v(x,0) & = & u_{0} (x)\quad \quad \quad in \quad \R.
\end{array} \right.
\end{eqnarray}
We now study the asymptotic development of the
solution of equation (\ref{gene1}).\\
The first goal of this paper is to analize the linear equation
\begin{equation}\label{genelin}
\left\{ 
\begin{array}{ll}
v_t + M v_t + Mv + Mv_x = 0, \quad x \in \R, \quad t>0;\\
\quad \; \;v(x,0) = u_0 (x),
\end{array}
\right.
\end{equation}
and to obtain the asymptotic expansion of solutions, complete when $m$ is 
an integer
and until the second term when $m$ is not integer respectively.\\
For the heat equation $u_t - \Delta u = 0$ in $\R^n$, 
this was done by J. Duoandixoetxea and E. Zuazua in \cite{Duozuazua}. 
Indeed, it was shown that if $u(x,t)=(G(t)\ast u_0)(x)$ is the solution 
of the heat equation whit initial data $u_0\in \Le^1(\Rn)$ and 
$G(x,t)=(4\pi t)^{-\frac{n}{2}} e^{-\frac{|x|^2}{4t}}$ 
($G$ is called the heat kernel),
with $u_0 \in \Le^1 (1+|x|^k)$ such that $|x|^{k+1}u_0 (x) \in \Le^p (\Rn)$, 
$1\leq p \leq q\leq \infty$, then 
\begin{equation}\label{duozua}
\nor G(t)\ast u_0 - \sum\limits_{|\alpha|\leq k}\frac{(-1)^{|\alpha|}}{\alpha !}
\biggl( \int x^{\alpha} u_0 (x)dx \biggr) D^{\alpha} G(t)\nor_q 
\leq Ct^{-(\frac{k+1}{2})- \frac{n}{2}(\frac{1}{p}-\frac{1}{q})} \,
\big\|\,|x|^{k+1} u_0 \big\|_p,
\end{equation}
for all $t>0$.\\
This shows that, for the solutions of the heat equation, a complete asymptotic
expansion may be obtained by means of the moments of the initial data and using
the derivatives of the Gaussian heat kernel as reference profiles. \\
In this work we show that for the solution of 
the linearized equation (\ref{genelin}), besides the terms
$$\sum\limits_{\alpha=0}^{k}\frac{(-1)^{\alpha}}{\alpha !}
\biggl( \int x^{\alpha} v_0 (x)dx \biggr) \partial_x^{\alpha} G_m(t)\;,$$ 
which correspond to the asymptotic expansion of the generalized linear heat 
equation
$$
v_t -Mv =0,
$$  
whit initial data $v_0 \in \Le^1(\R)$ and where
$G_m(x,t) = (1/(2\pi))\int_{\R} \text{e}^{ix.\xi-t|\xi|^m}d\xi$,  
other terms due to the dispersive effects appear in its asymptotic expansion.\\
We shall denote by $S(t)v_0$ the solution 
to the equation (\ref{genelin}), the function $K_m$ is defined through its  
Fourier transform $\fo{K_m}(\xi)= 1/(1+|\xi|^m)$, and 
$ K_m^j = \underbrace{K_m \ast K_m \ast ...\ast K_m}_{\rm j-veces}$.\\
Denoted
$${\cal{M}}_{\alpha}(v_0)=\frac{(-1)^{|\alpha|}}{\alpha !}
\biggl( \int x^{\alpha} v_0(x) dx \biggr)\;, $$
The following Theorems holds:
\begin{theo}\label{theo1}
Let $N \in \N$ and $m \in {\Z}^{+}$. Then there exist a constant $C=C(N)>0$ 
such that
\begin{multline*}
\nor {S}(t)v_0 \,- \!
\sum_{\alpha =0}^N  {\cal{M}}_{\alpha}(v_0) \partial_x^{\alpha}G_m(t) 
\,-\,\sum_{r=0}^{N}\frac{t^r}{r!}
\sum_{j=0}^{\maxi \frac{N}{2} \maxd}\frac{t^{j}}{j!}
\sum_{\scriptstyle 0\leq \alpha \leq N-r-mj 
\atop \scriptstyle (r,j)\not= 
(0,0)}^{*} 
\!\!\!\!\!{\cal{M}}_{\alpha}(v_0)
\big( \partial_x M \big)^r M^{2j} \partial_x^{\alpha}G_m(t) 
\nor_2 \\
\leq C \Bigg[ 
t^{- (\maxi \frac{N}{2} \maxd + 1) -\frac{1}{2m}} \big\| v_0 \big\|_1  
\biggl( \!\sum_{r=0}^{N} \frac{t^{-\frac{r}{m}}}{r!}\biggr)
+t^{-\frac{N+1}{m}-\frac{1}{2m}} \big\| v_0 \big\|_1
+t^{-\frac{N+1}{m}-\frac{1}{2m}}\sum_{r=1}^{N+1}\!
\big\| |x|^r v_0 \big\|_1  
+ e^{-\frac{t}{2}} \big\| v_0 \big\|_2    \\
+ e^{-\frac{t}{2}} t^{-\frac{1}{2m}}
\biggl( \sum_{r=0}^{N}\frac{t^{-\frac{r}{m}}}{r!}\biggr)
\biggl( \sum_{j=0}^{\maxi \frac{N}{2} \maxd}
\frac{t^{-j}}{j!}\biggr)  \big\| v_0 \big\|_1  
 +\max_{ 0\leq \alpha \leq N-1 } 
 \biggl\{  \big| {\cal{M}}_{\alpha}(v_0) \big|  \biggr\}
t^{-\frac{1+N}{m} - \frac{1}{2m}}
\big\| |x| K_{m} \big\|_1
\Bigg].
\end{multline*}
for all $v_0\in \Le^2(\R) \cap \Le^1(\R;1+|x|^{N+1})$\;.
\end{theo}
When $m$ is not integer, we still have the asymptotic expansion
until the second term of $S(t)v_0$, indeed:
\begin{theo}\label{theo2}
Let $N \in \N$ and $m = n+\delta$ for $n \in {\Z}^{+}$,$n>1$ $0 < \delta < 1$.
Then there exist a constant $C=C(N)>0$ such that
\begin{multline*}
\nor {S}(t)v_0  -
\sum_{\alpha =0}^{1} {\cal{M}}_{\alpha}(v_0)
\partial_x^{\alpha}G_m(t) -t {\cal{M}}_{0}(v_0) 
\big( \partial_x M \big) G_m(t) \nor_2  \\
\leq C \Bigg[
t^{-\frac{2}{m} - \frac{1}{2m}}  \big\|v_0 \big\|_1 
+t^{-\frac{2}{m} - \frac{1}{2m}}\sum_{r=1}^{2} \big\| |x|^r v_0 \big\|_1
+\;t^{-1-\frac{1}{2m}}
\biggl( \sum_{r=0}^{1} \frac{t^{-\frac{r}{m}}}{r!}\biggr) \big\| v_0 \big\|_1
+e^{-\frac{t}{2}} \big\| v_0 \big\|_2 \\
+ e^{-\frac{t}{2}}t^{-\frac{1}{2m}}
\biggl( \sum_{r=0}^{1} \frac{t^{-\frac{r}{m}}}{r!}\biggr) \big\| v_0 \big\|_1
+\;{\cal{M}}_{0}(v_0)\;t^{-\frac{2}{m} - \frac{1}{2m}}
\big\| |x| K_{m} \big\|_1
\Bigg]
\end{multline*}
for all $v_0\in \Le^2(\R) \cap \Le^1(\R;1+|x|^2)$\;.
\end{theo}
\begin{remark}
\item{(i) 
In the Theorem \ref{theo1}, for a real number 
$s>0$, $\maxi s \maxd = \max \{ m \in\N\; ; m \leq s\}$ denotes its  
integer part and
$\sum\limits_{|\alpha|\leq N-m-2j}^{*}$ means simply that the couples 
$(r,j)$ such that $N-r-mj < 0$ are not being considered in the sum. }
\item{(ii)
Our result generalizes that obtained in \cite{Karch2} and \cite{razu} 
where the case $m=2$ in the equation (\ref{genelin}) is studied.}
\item{(iii)
In the Theorems \ref{theo1}, when $m$ is a integer, we 
see two different terms in the expansion of $S(t)v_0$. The first one,  
$\sum\limits_{\alpha =0}^N  {\cal{M}}_{\alpha}(v_0) \partial_x^{\alpha}G_m(t)$ 
appears also 
in the asymptotic expansion of the generalized heat equation. The second one is 
due to the dispersive phenomena.}
\item{(iv)
In the Theorems \ref{theo2}, when $m$ is not a integer, we 
see two different terms in the expansion of $S(t)v_0$. The first one,  
$\sum\limits_{\alpha =0}^1  {\cal{M}}_{\alpha}(v_0) \partial_x^{\alpha}G_m(t)$ 
corresponding
to the two first term in the asymptotic expansion of the generalized heat equation.
The second one is $t {\cal{M}}_{0}(v_0) \big( \partial_x M \big) G_m(t)$
due to the dispersive phenomena caused by the operator $\partial_x M$.}
\end{remark}
We now study the nonlinear problem. That is, we consider the problem 
(\ref{gene1}), whit the following basic assumption: $q>m$, $m>2$,and that the 
solutions of (\ref{gene1}) satisfying the following decay estimates
\begin{gather} \label{optimos}
\big\|v(t) \big\|_2 \leq C(1+t)^{-\frac{1}{2m}}\;,\;\;\;
\big\|v_x(t) \big\|_2 \leq C(1+t)^{-\frac{1}{2m}-\frac{1}{m} }\;,\;\;\;
\big\|v(t) \big\|_{\infty} \leq C(1+t)^{-\frac{1}{m}},
\end{gather}
for all $t>0$, the numbers $C$ are independent of $t$.\\
We compute the second term of the asymptotic expansion of their solutions of the 
problem (\ref{gene1}) when $t\rightarrow \infty$.\\
When $m=2$, in \cite{Karch2} it was shown the influence of the nonlinear term
in the asymptotic expansion of their solutions. A result analogue when $m>2$
is the following Theorem: 
\begin{theo}\label{theo3}
Let $p\in [1,\infty)$, $q>m$, and $m>2$, and denote ${\cal{M}}=\int v_0(x)dx$.
Suppose that $u$ is a solution to
(\ref{gene1}), satisfying the decay estimates (\ref{optimos}), with 
$v_0 \in \Le^1 (\R)\cap H^2(\R)$, then
\begin{align*}
i)&\quad \text{For}\quad m <q< m+1,\\
& t^{(q-\frac{1}{p})\frac{1}{m}+\frac{1}{m}}
\nor v(t) - S(t)v_0 + 
\int_0^t \partial_x G_m(t-\tau) * ({\cal{M}} G_m (\tau))^q d\tau \nor_p 
\rightarrow 0\;,\quad \text{when}\quad t \rightarrow \infty\;.\\
ii)&\quad \text{For}\quad q = m+1,\\
& \frac{t^{(1-\frac{1}{p})\frac{1}{m}+\frac{1}{m}}}{\log t}
\nor v(t) - S(t)v_0 + 
\log t \bigg(\!\! \int \!\! \Big( \! {\cal{M}}G_m(t) \!\Big)^{m+1} \!\!(x,1)dx \bigg) 
\partial_x G_m(t) \nor_p 
\rightarrow 0,\quad \text{when}\;\; t \rightarrow \infty. \\
iii)&\quad \text{For}\quad q> m+1,\\
& t^{(1-\frac{1}{p})\frac{1}{m}+\frac{1}{m}}
\nor v(t) - S(t)v_0 + 
\bigg( \int_0^{\infty} \int_{\R} v^q (y,\tau)dyd\tau \bigg)
\partial_x G_m(t) \nor_p 
\rightarrow 0\;, \quad \text{when}\quad t \rightarrow \infty\;. \\
\end{align*}
\end{theo}
Of course, this result is to be complemented with Theorem \ref{theo1} and 
\ref{theo2} 
that providesa complete expansion of the linear component of $S(t)v_0$ to 
obtain
a complete description of the first and second terms of $v$.\\
At this respect, we recall that, when $n=1$ and $m=2$, G. Karch in \cite{Karch4} 
studied the equation (\ref{gene1})
and obtained the second order term for the cases of $2<q<3$, $q=3$, $q>3$. In this 
case 
$\int \!\! \big( \! {\cal{M}}G_2(t) \!\big)^{3} (x,1)dx = 
{\cal{M}}^3/(4\pi \sqrt{3})$.
\begin{remark}
The condition (\ref{optimos}) is not particularly restrictive, and is imposed for 
brevity sake. Indeed, see for example the work \cite{BP}. Moreover, for the 
completeness of our exposition, we construct solutions to (\ref{gene1}) 
satisfying (\ref{optimos}) provides $v_0$ is small in $W^{2,1}(\R)$ (see section 6).
\end{remark}
%
%
\section{Linearized equation}
In this section, we consider the family linearized (\ref{genelin}) of 
dispersive equations under the effect of dissipation.\\
Using the Fourier transform we obtain that each solution to (\ref{genelin}) 
with $u_0 \in {\cal{S}'}$ has the form
\begin{gather}
v(x,t)= S(t)v_0 (x) = 
\frac{1}{2\pi} \int_{\R} \text{e}^{-t\Phi{\xi}+ix.\xi} \fo{v_0}{\xi} d\xi
\end{gather}
whit the phase function 
$$
\Phi(\xi) = \frac{|\xi|^m - i\xi|\xi|^m}{1 + |\xi|^m}
$$
Our main goal will be to find the complete asymptotic expansion of solutions
$S(t)v_0 (x)$ as $t$ approaches infinity, with this aim, first 
we obtain in the following subsection, the complete
asymptotic expansion of solutions, as time goes to infinity, of the 
generalized heat equation.
%
\subsection{Generalized heat equation. 
Complete asymptotic expansion.}
%
The complete asymptotic expansion for the solution of heat equation were 
obtained by J. Duoandikoetxea and E. Zuazua in \cite{Duozuazua}, in similar 
form, we have the complete asymptotic expansion for a generalization of the 
heat equation, our generalization takes the form
\begin{equation}\label{calge}
\left\{
\begin{array}{ll}
u_t + M u = 0, \quad x \in \R, \quad t>0;\\
\quad \; \;u(x,0) = u_0 (x),
\end{array}
\right.
\end{equation}
Where the operator $M$ is given by (\ref{opM}), with $m\geq 1$.\\
If $u_0 \in \Le^1 (\R)$, using the Fourier transform in the variable $x$,
we obtain that the solution to (\ref{calge}) is a convolution product:\;
$u(x,t) = G_m(t)\ast u_0 (x)$, where
\begin{equation}
G_m (x,t) = \frac{1}{2\pi} \int_{\R} e^{-|\xi|^{m} t + ix \cdot \xi} d\xi,
\end{equation}
is the heat kernel generalized.\\
If $t>0$ and $j$ is a nonnegative integer, then using the Euler integral
it can be easily seen that
$$
\int_{\R} \xi^{2j} e^{-2t\xi^{m}}d\xi = C(m,j) t^{-\frac{2j}{m} - \frac{1}{m}},
\quad\text{where}\quad C(m,j)= \frac{\Gamma(\frac{2j+1}{m})}{m}.
$$
Hence, it follows that
\begin{equation}
\big\| \partial_x^j G_m (t) \big\|_2 = C(m,j) t^{-\frac{j}{m} - \frac{1}{2m}}
\end{equation}
In general we have the following result:
\begin{equation}\label{nucleo}
\big\| \partial_x^j G_m (t) \big\|_p 
\leq C t^{-\frac{j}{m} - \frac{1}{m}(1-\frac{1}{p})}, 
\quad\forall \;1 \leq p \leq \infty, \quad m \geq 1.
\end{equation}
For the proof of (\ref{nucleo}) we estimate 
$\big\| \partial^j G_m (t) \big\|_{\infty}$ and 
$\big\|\partial^j G_m (t)\big\|_1$, then 
interpolation. The case $p=1$, use the following formula (ver \cite{Karch1} )
$$
\big\| w \big\|_1 \leq 
C \big\| w \big\|_2^{1/2} \big\| \partial_{x}\fo{w} \big\|_2^{1/2}.
$$
We have the complete asymptotic expansion for the solutions of (\ref{calge}) 
by means of the moments of the initial data and using the derivatives of the 
generalized heat kernel as reference profiles. Indeed, if denote
the moments of the initial data by 
$$
{\cal{M}}_j (u_0)=  \frac{(-1)^j}{j!} \int_{\R} x^j u_0(x)dx,
$$
we have the following Theorem:
\begin{theo}\label{calorg}
Assume that $u$ is the solution to (\ref{calge}), with $u_0 \in \Le^1(\R)$.
Let $j,N$ nonnegatives integer. Then, there exist a constant $C>0$ such that
\begin{equation}
\nor G_m(t)\ast u_0 - \sum_{j=0}^N {\cal{M}}_j(u_0)\partial_x^j G_m(t) \nor_p
\leq C t^{-\frac{1}{m}(1-\frac{1}{p})-\frac{N+1}{m}} 
\big\| |x|^{N+1} u_0 \big\|_1 , \quad \forall \, p\in[1,\infty]\,,
\end{equation}
for all $u_0 \in \Le((1+|x|^N)dx,\R) \cap \Le^1 (|x|^{N+1}dx,\R)$
\end{theo}
\proof
By \cite{Duozuazua} $u_0$ can be decomposed as
$$
u_0 = \sum_{j=0}^N {\cal{M}}_j(u_0) \partial_x^j \delta + \partial_x^{N+1} F_{N+1},
$$
where $F_{N+1} \in \Le^1(\R)$, such that 
$\big\| F_{N+1} \big\|_1 \leq \big\| |x|^{N+1}u_0 \big\|_1$.\\
Hence, taking the convolution of $u_0$ with $G_m(t)$, we obtained
\begin{multline*}
\nor G_m(t)\ast u_0 - \sum_{j=0}^N {\cal{M}}_j(u_0)\partial_x^j G_m(t)\nor_p 
= \nor \partial^{N+1} F_{N+1} \ast G_m(t) \nor_p  
\leq \big\| F_{N+1} \big\|_1 \big\| \partial^{N+1} G_m(t)\big\|_p \\
\leq C t^{-\frac{N+1}{m}-\frac{1}{m}(1-\frac{1}{p})} 
\big\| |x|^{N+1}u_0\big\|_1 . \hspace{1.5cm}{\BO}
\end{multline*}
The following result gives the different representations of the operator $M$
\begin{lem}(\cite{CBona})
Let $M$ the operator defined by $\fo{Mu}(\xi)=|\xi|^m \fo{u}(\xi)$, for 
$u \in H^r$, where $0 < m \leq r$, $r$ is a positive integer. It follows that
\begin{enumerate}
\item let $m=2n$ for $n \in Z^{+}$, then
$$
Mu(x)= (-1)^n \partial_{x}^{2n} u(x);
$$
\item let $m=2n+1$ for $n \in Z^{+}$, then
$$
Mu(x)=(-1)^n \sqrt{\frac{2}{\pi}} \int \frac{\partial_y^{2n+1}u(y)}{x-y} dy;
$$
\item let $m=2n+\delta$ for $n \in Z^{+}$, $0 < \delta < 1$, then
$$
Mu(x)=(-1)^n \sqrt{\frac{\pi}{2}} 
\Big ( cos\big(\frac{\delta\pi}{2}\big)\Gamma(1-\delta)\Big)^{-1} 
\int sign(x-y)\frac{\partial_y^{2n+1}u(y)}{|x-y|^{\delta}} dy,
$$
where $\Gamma$ denote the gamma function;
\item let $m=2n+1+\delta$ for $n \in Z^{+}$, $0 < \delta < 1$, then
$$
Mu(x)=(-1)^n \sqrt{\frac{\pi}{2}}
\Big ( sen\big(\frac{\delta\pi}{2}\big)\Gamma(1-\delta)\Big)^{-1} 
\int \frac{\partial_y^{2n+1}u(y)}{|x-y|^{\delta}} dy.
$$
\end{enumerate}
\end{lem}
%
\subsection{Complete asymptotic expansion: linear case}
%
In this section we obtain preliminary result on the asymptotic development of
solutions $S(t)v_0$ of (\ref{genelin}). In this developmente, appear terms of 
convolution
$K_m^{r+j}\ast (\partial_x M)^r M^{2j} \partial_x^{\alpha}G_m(t)$, then in 
the following section, to prove Theorems \ref{theo1} and \ref{theo2} we replace
them by the first term of their asymptotic expansion, i.e. by 
$(\partial_x M)^r M^{2j} \partial_x^{\alpha}G_m(t)$.\\
The following preliminary Theorems holds:
\begin{theo}\label{theoa}
Let $N \in \N$ and $m \in {\Z}^{+}$.
Then, there exist a constant
$ C = C(N) > 0$ such that
\begin{align*}
\hspace{-1.5cm}\nor S(t)v_0  & -
\sum_{r=0}^{N}\frac{t^r}{r!}\sum_{j=0}^{\maxi \frac{N}{2} \maxd}
\frac{t^{j}}{j!}
\sum_{\scriptstyle 0\leq\alpha\leq N-r-mj}^{*}
{\cal{M}}_{\alpha}(v_0)
 K_m^{r+j}\ast (\partial_x M)^r M^{2j} \partial_x^{\alpha}G_m(t)\nor_2  \\
& \leq Ct^{-(N+1)\frac{1}{m} - \frac{1}{2m}}  \big\| v_0 \big\|_1 
+\;Ct^{-(\maxi \frac{N}{2} \maxd +1)-\frac{1}{2m}}
\biggl( \sum_{r=0}^{N} \frac{t^{-\frac{r}{m}}}{r!}\biggr)\big\| v_0 \big\|_1 \\
&\quad\quad \quad  + Ct^{-(N+1)\frac{1}{m}-\frac{1}{2m}}
\sum_{r=0}^{N}
\sum_{\scriptstyle 0\leq j \leq \maxi \frac{N}{2} \maxd 
\atop \scriptstyle 0\leq j \leq \maxi \frac{N-r}{m} \maxd}
\Big\| |x|^{N+1-r-jm}v_0 \Big\|_1 \\
&\quad \quad\quad\quad + Ce^{-\frac{t}{2}}\big\| v_0 \big\|_2 +
Ce^{-\frac{t}{2}} t^{-\frac{1}{2m}}
\biggl( \sum_{r=0}^{N}\frac{t^{-\frac{r}{m}}}{r!}\biggr)
\biggl( \sum_{j=0}^{\maxi \frac{N}{2} \maxd} \frac{t^{-j}}{j!}\biggr)
\big\| v_0 \big\|_1 \:
\end{align*}
for all $v_0\in \Le^2(\R) \cap \Le^1(\Rn;1+|x|^{N+1})$\;.
\end{theo}
When $m$ is not integer, we still have the asymptotic expansion
until the second term of $S(t)v_0$ indeed, we have:
%
\begin{theo}\label{theob}
Let $N \in \N$ and $m =n + \delta$, for $n\in {\Z}^{+}$, $0<\delta<1$.
Then, there exis a constant
$ C = C(N) > 0$ such that
\begin{align*}
\nor S(t)v_0  -
\sum_{\alpha =0}^{1} {\cal{M}}_{\alpha}(v_0) &
\partial_x^{\alpha}G_m(t) -tMK_m\ast\partial_x M \nor_2  \leq 
Ct^{-\frac{2}{m} - \frac{1}{2m}}  \big\| v_0 \big\|_1 \\
&+Ct^{-\frac{2}{m} - \frac{1}{2m}}\sum_{r=1}^{2}\big\| |x|^r v_0 \big\|_1
+\;Ct^{-1-\frac{1}{2m}}
\biggl( \sum_{r=0}^{1} \frac{t^{-\frac{r}{m}}}{r!}\biggr) \big\| v_0 \big\|_1 \\
&\quad \quad +Ce^{-\frac{t}{2}}\big\| v_0\big \|_2 +
Ce^{-\frac{t}{2}}t^{-\frac{1}{2m}}
\biggl( \sum_{r=0}^{1} \frac{t^{-\frac{r}{m}}}{r!}\biggr) \big\| v_0 \big\|_1,
\end{align*}
for all $v_0\in \Le^2(\R) \cap \Le^1(\Rn;1+|x|^2)$\;.
\end{theo}
For the proof of Theorems \ref{theoa} and \ref{theob}, we will use the following
Lemmas. We consider $\varphi\in {\cal{C}}_{0}^{\infty}(\R)$, a cut-off function
such that:
\begin{align*}
|\varphi(\xi)|&\leq 1 ,\quad \forall\xi\in \Rn \;;\\
\varphi(\xi)& = 1 ,\quad  \:\text{if} \quad |\xi|\leq 1 \;;\\
\varphi(\xi)& = 0 , \quad \:\text{if} \quad |\xi| > 2. \\
\end{align*}
Then, we define
\begin{equation}\label{sphi}
S_{\varphi}(x,t)= 
\frac{1}{2\pi}\int e^{-t \Phi_0(\xi) + ix\cdot \xi} \varphi(\xi)d\xi.
\end{equation}
Where
$$
\Phi_0 (\xi) = \frac{|\xi|^m}{1 + |\xi|^m}, \quad m \geq 1.
$$
Notice that if define
$$
K_m^j = \underbrace{K_m \ast K_m \ast ...\ast K_m}_{\rm j-veces},
$$
where $\fo{K}_m (\xi) = 1/(1 + |\xi|^m)$, then
\begin{gather}\label{funK}
\big[ \fo{K}_m (\xi) \big]^j = \fo{K_m^j}(\xi) = 
\frac{1}{\big[ 1 + |\xi|^m \big]^j}.
\end{gather}
Hence
\begin{align*}
\frac{1}{2\pi}
\int e^{-t \Phi_0(\xi) + ix\cdot \xi} &\;\fo{u}_0(\xi)d\xi =
\frac{1}{2\pi}
\int e^{-t|\xi|^m+ ix\cdot \xi}\; e^{\frac{t|\xi|^{2m}}{1+|\xi|^m}} \;\fo{u}_0(\xi)d\xi\\
&=\sum_{j=0}^{\infty}\frac{1}{2\pi} \int e^{-t|\xi|^m + ix\cdot \xi}
\bigg(t |\xi|^{2m}\;\fo{K_m}{\xi}\bigg)^j \;\frac{1}{j!}\;\fo{u}_0(\xi) d\xi\\
&= \sum_{j=0}^{\infty}\frac{t^j}{j!}\big( K_m^j \ast M^{2j} G_{m}\ast u_0 \big)(x).
\end{align*}
Then, we have the following results:
\begin{lem}\label{lemaa} 
Denote $K_m^j = \underbrace{K_m \ast K_m \ast ...\ast K_m}_{\rm j-veces}$, 
$m \geq 1$. 
Let $N \in \N$. There exist a constant $C=C(N,m)$ such that 
\begin{multline}\label{lem0}
\nor S_{\varphi}(t)\ast u_0 - \sum_{j=0}^{N} \frac{t^j}{j!}
K_m^j \ast M^{2j} G_{m}\ast u_0 \nor_2 \leq 
C t^{ -(N+1) - \frac{1}{2m} } \big\| u_0 \big\|_1  \\
+ C e^{-t/2} t^{-\frac{1}{2m}}
\bigg( \sum_{j=0}^{N}\frac{t^{-j}}{j!} \bigg) \big\| u_0 \big\|_1\,,
\end{multline}
for all $t>0$ and $u_0 \in \Le^1 (\R)$.
\end{lem}
\proof
We decompose $G_m(x,t)$ using the cut-off function $\varphi(\xi)$
$$
G_m (x,t)= G_{m\varphi} (x,t)+ G_{m(1-\varphi)}(x,t).
$$
Then
\begin{multline}\label{lem1}
\nor S_{\varphi}(t)\ast u_0 - \sum_{j=0}^{N} \frac{t^j}{j!}
K_m^j \ast M^{2j} G_{m}\ast u_0 \nor_2 \leq
 \nor \sum_{j=0}^{N} \frac{t^j}{j!}K_m^j \ast M^{2j} G_{m(1-\varphi)}\ast u_0) \nor_2 \\
+ \nor S_{\varphi}(t)\ast u_0 - \sum_{j=0}^{N} \frac{t^j}{j!}
K_m^j \ast M^{2j} G_{m\varphi}\ast u_0 \nor_2 = I_1(t) + I_2(t).
\end{multline}
We estimate $I_i(t)$, $i=1,2$, separately.\\
{\it The term $I_2(t)$}: Let
$$
w(x,t) = S_{\phi}(t) - \sum_{j=0}^{N} \frac{t^j}{j!}K_m^j \ast M^{2j} G_{m\varphi}.
$$
Then, its Fourier transform is
\begin{align}\label{trans}
\fo{w}(\xi,t)&= e^{-\frac{t|\xi|^m}{1+|\xi|^m}} \varphi(\xi)-
\sum_{j=0}^{N} \frac{t^j}{j!} 
\bigg( \frac{|\xi|^{2m}}{1+|\xi|^m} \bigg)^j e^{-t|\xi|^m} \varphi(\xi) \\
&= e^{-t|\xi|^m} \varphi(\xi) \Bigg[ e^{\frac{t|\xi|^{2m}}{1+|\xi|^m}} -
\sum_{j=0}^{N} \frac{t^j}{j!} 
\bigg( \frac{|\xi|^{2m}}{1+|\xi|^m} \bigg)^j \Bigg]\;.\nonumber
\end{align}
Hence, using the Taylor expansion of the exponential function, we have
$$
e^x -\sum_{j=0}^N \frac{x^j}{j!}\leq \frac{x^{N+1}}{(N+1)!} e^x.
$$ 
Using this inequality in (\ref{trans}) we obtain
\begin{align}\label{wfou}
\big| \fo{w}(\xi,t) \big| &\leq 
e^{-t|\xi|^m} \varphi(\xi) \bigg( t\frac{|\xi|^{2m}}{1+|\xi|^m}\bigg)^{N+1}
e^{\frac{t|\xi|^{2m}}{1+|\xi|^m}}  \\
&\leq C_N \;\varphi(\xi)\;e^{-\frac{t|\xi|^m}{1+|\xi|^m}}\;
\bigg( t\frac{|\xi|^{2m}}{1+|\xi|^m}\bigg)^{N+1}.\nonumber
\end{align}
By Plancherel's formula we see that
\begin{align}\label{wplan1}
\big\| w(t)\big\|_2^2 &= \int \big| \fo{w} \big|^2 d\xi 
\leq C_N \int |\varphi(\xi)|^2 e^{-\frac{2t|\xi|^m}{1+|\xi|^m}}
\bigg( \frac{t|\xi|^{2m}}{1+|\xi|^m}\bigg)^{2(N+1)} d\xi \\
&\leq C t^{2(N+1)} \int_{|\xi|\leq 1}e^{-\frac{2t|\xi|^m}{1+|\xi|^m}} 
\;\frac{|\xi|^{4m(N+1)}}{1+|\xi|^m} d\xi \nonumber \\
\leq& C t^{2(N+1)} \int e^{-t|\xi|^m} \big| \xi \big|^{4m(N+1)}\; d\xi,
\:(m\geq 1,\:|\xi|\leq 1 \Rightarrow \frac{2|\xi|^m}{1+|\xi|^m}\geq |\xi|^m)\nonumber \\
&\leq C t^{2(N+1)} t^{-\frac{4m(N+1)}{m} - \frac{1}{m}}=
C t^{ -2(N+1) - \frac{1}{m} }.\nonumber
\end{align}
Then, from (\ref{wplan1}) we have
\begin{equation}\label{dei2}
I_2(t)\leq \nor S_{\varphi}\ast u_0 - \sum_{j=0}^{N} \frac{t^j}{j!}
K_m^j \ast M^{2j} G_{m\varphi} \nor_2 \; \big\| u_0 \big\|_1 
\leq C t^{ (N+1) - \frac{1}{2m} } \big\| u_0 \big\|_1.
\end{equation}
{\it The term $I_1(t)$}:\\
Given that $\big\| K_m^j \big\|_1 = 1$, then we have
\begin{equation}\label{wplan2}
I_1(t)\leq \sum_{j=0}^{N}\frac{t^j}{j!} 
\big\| K_m^j \big\|_1 \; \big\| M^{2j} G_{m(1-\varphi)} \big\|_2 \; \big\| u_0 \big\|_1
\leq  \sum_{j=0}^{N}\frac{t^j}{j!} \big\| M^{2j}
 G_{m(1-\varphi)} \big\|_2\big\|\; u_0 \big\|_1.
\end{equation}
Now, by Plancherel's formula and using that
$|\xi| \geq 1$ implies $2|\xi|^m \geq 1 + |\xi|^m$, we have
\begin{align}\label{wplan3}
\big\| M^{2j} G_{m(1-\varphi)} \big\|_2^2 &= 
\int \big| 1-\varphi(\xi) \big|^2 \big| \xi \big |^{4jm} e^{-2t|\xi|^m}d\xi 
\leq \int_{|\xi|\geq 1/2} \big| \xi \big|^{4mj}e^{-2t|\xi|^m}d\xi \\
&\leq e^{-t} \int e^{-t|\xi|^m}\;|\xi|^{4mj} d\xi 
\leq C e^{-t}t^{-4j-\frac{1}{m}}.\nonumber
\end{align}
Then, using (\ref{wplan3}) in (\ref{wplan2}) we have
\begin{equation}\label{wplan4}
I_1(t)\leq C e^{-t/2} t^{-\frac{1}{2m}}\;
\sum_{j=0}^{N}\frac{t^{-j}}{j!} \big\| u_0 \big\|_1.
\end{equation}
Then, from (\ref{wplan4}) and (\ref{dei2}) we have (\ref{lem0}).
\hfill$\BO$
%
\begin{lem}\label{lemab}
Let $N \in \N$. Then there exist a constant positive $C=C(N,m)$ such that
\begin{equation}
\nor S(t)u_0 - 
\sum_{j=0}^{N} \frac{t^j}{j!} K_m^j \ast \big( \partial_x M \big)^j S_{\varphi}\ast u_0 \nor_2
\leq C \Bigl[ t^{-\frac{N+1}{m} - \frac{1}{2m}} \big\| u_0 \big\|_1 + 
e^{-t/2}\big\| u_0 \big\|_2 \Bigr] ,
\end{equation}
for all $t>0$ and $u_0 \in \Le^2 (\R) \cap \Le^1 (\R)$
\end{lem}
\proof
We decompose $(t)$ using the cut-off function $\varphi$, we have 
\begin{multline}\label{jota}
S(t)u_0 - \sum_{j=0}^{N} \frac{t^j}{j!} K_m^j 
\ast \big( \partial_x M \big)^j S_{\varphi}\ast u_0 =
\frac{1}{2\pi} \int e^{t\widetilde{\Phi} + ix\cdot\xi} \;\fo{u_0}(\xi) 
\big( 1-\varphi(\xi)\big) \;d\xi \\
+ \frac{1}{2\pi} \int e^{t\widetilde{\Phi} + ix\cdot\xi} \;\fo{u_0}(\xi)\varphi(\xi)d\xi -
\sum_{j=0}^{N} \frac{t^j}{j!} K_m^j \ast 
\big( \partial_x M \big)^j S_{\varphi}\ast u_0 (x)
= J_1(t)+J_2(t).
\end{multline}
Here $J_2(t)$ represents the difference of the last two terms in (\ref{jota}).
We estimate $J_i(t)$, $i=1,2$, separately.\\
{\it The term $J_2(t)$}:\\
The Fourier transform of $J_2(x,t)$ is
\begin{align*}
\fo{J_2}(\xi,t) &=  e^{-t\widetilde{\Phi}} \varphi(\xi)\fo{u_0}(\xi) - 
\sum_{j=0}^{N} \frac{t^j}{j!} \bigg( \frac{-i\xi |\xi|^m}{1+|\xi|^m} \bigg)^j 
e^{-\frac{t|\xi|^m}{1+|\xi|^m}} \varphi(\xi) \fo{u_0}(\xi)\\
&=  e^{-\frac{t|\xi|^m}{1+|\xi|^m}}\varphi(\xi)\fo{u_0}(\xi) 
\Bigg[ 
e^{-\frac{it\xi |\xi|^m}{1+|\xi|^m}} - 
\sum_{j=0}^{N} \frac{t^j}{j!} \bigg( -\frac{i\xi|\xi|^m}{1+|\xi|^m}\bigg)^j 
\Bigg].
\end{align*}
The Taylor expansion of $e^{ix}$ implies that 
$\bigl|e^{ix} -\sum_{j=0}^N \frac{(ix)^j}{j!}\bigr| \leq \frac{x^{N+1}}{(N+1)!}$.
Then
\begin{equation}\label{lem2}
\big| \fo{J_2}(\xi,t) \big| \leq C e^{-\frac{t|\xi|^m}{1+|\xi|^m}}\;
\big| \varphi(\xi) \big|
\big| \fo{u_0}(\xi) \big|
\bigg( \frac{t\xi|\xi|^m}{1+|\xi|^m}\bigg)^{N+1}
\end{equation}
By Plancherel's formula, by (\ref{lem2}) and observing that  
$|\xi|\leq 1$ implies $\frac{2|\xi|^m}{1+|\xi|^m} \geq |\xi|^m,\;(m\geq 1)$,
we have
\begin{align}\label{lem3}
\big\| J_2(t) \big\|_2^2 &= \int \big| \fo{J_2}(\xi) \big|^2 d\xi \leq C_N 
\int e^{-\frac{2t|\xi|^m}{1+|\xi|^m}} 
\big|\varphi(\xi) \big|^2 \big| \fo{u_0}(\xi) \big|^2
\big( t\xi |\xi|^m \big)^{2(N+1)} d\xi \\
&\leq C_N t^{2(N+1)} \big\| u_0 \big\|_1^2
\int_{|\xi|\leq 1}  e^{-\frac{2t|\xi|^m}{1+|\xi|^m}} |\xi|^{2(m+1)(N+1)} \nonumber \\
&\quad\quad\leq C t^{2(N+1)}\big\| u_0 \big\|_1^2 \int e^{-t|\xi|^m}
\big| \xi \big|^{2(m+1)(N+1)} d\xi 
\leq C t^{-\frac{2(N+1)}{m}- \frac{1}{m}} \big\| u_0 \big\|_1^2.\nonumber
\end{align}
{\it The term $J_1(t)$}:\\
The Fourier transforme of $J_1(x,t)$ is
\begin{equation*}
\fo{J_1}(\xi) = e^{-t \varphi(\xi)}\; \fo{u_0}(\xi)\; \big( 1- \varphi(\xi) \big).
\end{equation*}
Then, by Plancherel's formula we have
\begin{align}\label{lem4}
\big\| J_1(t)\big\|_2^2 &= \int \big| \fo{J_1}(\xi,t) \big|^2 d\xi = 
\int e^{\frac{-2t|\xi|^m}{1+|\xi|^m}} 
\big| \fo{u_0}(\xi) \big|^2 \big| 1-\varphi(\xi) \big|^2 d\xi \\
& \leq \int_{|\xi|\geq 1/2} 
e^{\frac{-2t|\xi|^m}{1+|\xi|^m}} \big| \fo{u_0}(\xi) \big|^2 d\xi. \nonumber
\end{align}
Now, we have that if $m \geq 1$, then
$$|\xi| \geq 1 \Rightarrow 2|\xi|^m \geq 1 + |\xi|^m $$
Returning to (\ref{lem4}), we have 
\begin{equation}\label{lem5}
\big\| J_1(t) \big\|_2^2 \leq \int e^{-t} \big| u_0(\xi) \big|^2 d\xi = 
e^{-t}\big\| u_0 \big\|_2^2.
\end{equation}
The Lemma \ref{lemaa} is consequence of (\ref{lem3}) and (\ref{lem5}).
$\hfill \BO$
%
%
\begin{lem}\label{lemac}
Let $ r\in \N$, $r\geq 1$. Then there exist a constant
$C=C(N,r) > 0$ such that
\begin{multline*}
\nor K_m^r\ast \big( \partial_x M \big)^r
S_{\varphi}(t) \ast v_0 -  
\sum_{j=0}^{N} \frac{t^j}{j!}
K_m^{r+j} \ast \big(\partial_x M \big)^r
M^{2j} G_m(t)\ast v_0  \nor_2 \\
\leq Ct^{-\frac{1}{m}(r(m+1)+m(N+1)) - \frac{1}{2m}} \big\| v_0 \big\|_1 + 
Ce^{-\frac{t}{2}} t^{-(\frac{m+1}{m})r - \frac{1}{2m}}
\biggl(\sum_{j=0}^{N} \frac{t^{-j}}{j!}\biggr) \big\| v_0 \big\|_1 
\end{multline*}
for all $t>0$ and $v_0 \in \Le^1 (\R)$.
\end{lem}
\proof
Given that $\big\| K_m^r \big\|_1 = 1$, we have
\begin{multline*}
\nor K_m^r\ast \big(\partial_x M \big)^r 
 S_{\varphi}(t)\ast v_0 -    
\sum_{j=0}^{N} \frac{t^j}{j!}
K_m^{r+j} \ast \big( \partial_x M \big)^r M^{2j} G_m(t)\ast v_0  \nor_2 \leq \\
\leq \nor \big( \partial_x M \big)^r S_{\varphi}(t)\ast v_0 - 
\sum_{j=0}^{N} \frac{t^j}{j!} 
 K_m^j\ast \big(\partial_x M \big)^r M^{2j} G_m(t)\ast v_0 \nor_2\;.
\end{multline*}
Hence, we only need to prove that 
\begin{multline}\label{2.3.1}
\nor \big( \partial_x M \big)^r S_{\varphi}(t)\ast v_0 - 
\sum_{j=0}^{N}\frac{t^j}{j!}
 K_m^j\ast \big( \partial_x M \big)^r M^{2j} G_m(t)\ast v_0 \nor_2 \leq  \\
\leq Ct^{-\frac{1}{m}(r(m+1)+m(N+1))- \frac{1}{2m}} \big\| v_0 \big\|_1 + 
Ce^{-\frac{t}{2}} t^{-(\frac{m+1}{m})r - \frac{1}{2m}}
\biggl(\sum_{j=0}^{N} \frac{t^{-j}}{j!}\biggr) \big\| v_0 \big\|_1  
\end{multline}
We decompose $G_m$ using the cut-off function $\varphi$:
$$
G_m(x,t)=G_{m\varphi}(x,t) + G_{m(1-\varphi)}(x,t).
$$
Then
\begin{align}\label{lem2.3.2}
\nor \big( \partial_x M \big)^r S_{\varphi}(t)\ast & v_0 -
\sum_{j=0}^{N}\frac{t^j}{j!}
 K_m^j\ast \big( \partial_x M \big)^r  M^{2j} G_m(t)\ast v_0 \nor_2 \leq \\
&\leq \nor \big(\partial_x M \big)^r S_{\varphi}(t)\ast v_0 
- \sum_{j=0}^{N}\frac{t^j}{j!}
 K_m^j\ast \big( \partial_x M \big)^r  M^{2j} G_{m\varphi}(t)\ast v_0 \nor_2 \nonumber\\
&\quad + \nor \sum_{j=0}^{N}\frac{t^j}{j!}
 K_m^j\ast \big( \partial_x M \big)^r M^{2j} G_{m(1-\varphi)}(t) \ast v_0 \nor_2 
=I_1 (t) + I_2 (t)\:.\nonumber
\end{align}
We estimate $I_1 (t)$ e $I_2 (t)$ separately.\\
{\it The term $I_1 (t)$:} We have
\begin{multline}\label{lem2.3.3}
\nor \big( \partial_x M \big)^r S_{\varphi}(t)\ast v_0 -
\sum_{j=0}^{N}\frac{t^j}{j!}
 K_m^j\ast \big( \partial_x M \big)^r M^{2j} G_{m\varphi}(t)\ast v_0 \nor_2 \leq \\
\leq \nor \big(\partial_x M \big)^r S_{\varphi}(t) - 
\sum_{j=0}^{N}\frac{t^j}{j!}
 K_m^j\ast \big( \partial_x M \big)^r M^{2j} G_{m\varphi}(t) \nor_2 \;
 \big\| v_0 \big\|_1\:.
\end{multline}
We now get, 
$g(x,t) =\big( \partial_x M \big)^r S_{\varphi}(t) - 
\sum\limits_{j=0}^{N}\frac{t^j}{j!}
 K_m^j\ast \big( \partial_x M \big)^r  M^{2j} G_{m\varphi}(t)$. 
Its Fourier transform is
\begin{eqnarray}\label{lem2.3.4}
\fo{g}(\xi,t)&=&
\big( i\xi |\xi|^m \big)^r\varphi (\xi)\Biggl(
e^{-\frac{t|\xi|^m}{1+|\xi|^m}}
\:-\:\sum_{j=0}^{N}\frac{t^j}{j!}\biggl( \fo{K_m}(\xi)\,\big| \xi \big|^{2m}\biggr)^j
\:e^{-t|\xi|^m}\Biggr)   \\ 
&=& \varphi (\xi) \big( i\xi |\xi|^m \big)^r  e^{-t|\xi|^m}
\Biggl[ e^ {\frac{t|\xi|^{2m}}{1+|\xi|^m}} \:-\:
\sum_{j=0}^{N}\frac{t^j}{j!}\biggl( \frac{|\xi|^{2m}}{1+|\xi|^m}
\biggr)^j\,\Biggr]\:.\nonumber
\end{eqnarray}
Now, using the Taylor expansion of the function
$e^x, x\geq 0$, we have
$e^x - \sum_{j=0}^N \frac{x^j}{j!} \leq \frac{x^{N+1}}{(N+1)!}e^x$.
Using this inequality in (\ref{lem2.3.4}) we have
\begin{eqnarray}\label{lem2.3.5}
\big| \fo{g}(\xi,t) \big|\, &\leq&C\,\varphi(\xi)\,\big| \xi \big|^{r(m+1)}
\frac{e^{-t|\xi|^m}}{(N+1)!}
\biggl(\frac{t|\xi|^{2m}}{1+|\xi|^m}\biggr)^{N+1}
e^{\frac{t|\xi|^{2m}}{1+|\xi|^m}}  \\
&=&C\,\varphi(\xi)\,\big| \xi \big|^{r(m+1)}
\frac{e^{-\frac{t|\xi|^m}{1+|\xi|^m}}}{(N+1)!}
\biggl(\frac{t|\xi|^{2m}}{1+|\xi|^m}\biggr)^{N+1}\;.\nonumber
\end{eqnarray}
>From (\ref{lem2.3.5}) and Plancherel's formula we obtain
\begin{align}\label{lem2.3.6}
\nor \big(\partial_x M \big)^r S_{\varphi}(t) & -
\sum_{j=0}^{N}\frac{t^j}{j!}
 K_m^j\ast \big(\partial_x M \big)^r M^{2j}  G_{m\varphi}(t) \nor_{2}^{2} = 
\int \big| \fo{g}(\xi,t) \big|^2 d\xi  \\
&\leq C_N\int \big| \varphi(\xi)\big|^2\,|\xi|^{2r(m+1)}
e^{-\frac{2t|\xi|^m}{1+|\xi|^m}}
\biggl(\frac{t|\xi|^{2m}}{1+|\xi|^m}\biggr)^{2(N+1)}d\xi \nonumber \\
&\quad =\; C_N t^{2(N+1)}\int_{|\xi|\leq 1}\,\big| \xi \big|^{2r(m+1) + 4m(N+1)}
e^{-t|\xi|^m}d\xi  \nonumber\\
& \quad\quad \leq C_N t^{2(N+1)}\times t^{-\frac{2}{m}[r(m+1) + 2m(N+1)] - \frac{1}{m}}
= Ct^{-\frac{2}{m}[r(m+1) + m(N+1)] - \frac{1}{m}}.\nonumber
\end{align}
In (\ref{lem2.3.6}) we use that
$|\xi|\leq 1$ implies $\frac{2|\xi|^m}{1+|\xi|^m} \geq |\xi|^m$.
Then, returning (\ref{lem2.3.3}), by (\ref{lem2.3.6}) it follows that
\begin{multline}\label{lem2.3.7}
\nor \big( \partial_x M \big)^r S_{\varphi}(t)\ast v_0  -
\sum_{j=1}^{N}\frac{t^j}{j!}
 K_m^j\ast \big(\partial_x M \big)^r M^{2j} G_{m\varphi}(t)\ast v_0 \nor_2 \leq  \\
\leq \;Ct^{-\frac{1}{m}[r(m+1) + m(N+1)] - \frac{1}{2m}} \big\| v_0 \big\|_1\;.
\end{multline}
{\it The term $I_2(t)$:}
Since $\big\|  K_m^j \big\|_1 = 1$, we have
\begin{equation}\label{lem2.3.8} 
\nor \sum_{j=0}^{N}\frac{t^j}{j!}
 K_m^j\ast \big(\partial_x M \big)^r M^{2j} G_{m(1-\varphi)}
\ast v_0 \nor_2 
\leq \;\sum_{j=0}^{N}\frac{t^j}{j!}
\big\| \big(\partial_x M \big)^r M^{2j} G_{m(1-\varphi)} \big\|_2 \big\| v_0 \big\|_1\;,
\end{equation}
and moreover if $m \geq 1$, 
$|\xi| \geq 1$ implies that $2t\big| \xi \big|^m \geq \big(t + t|\xi|^m \big)$.
Then
\begin{align*}
\Big\| \big( \partial_x M \big)^r M^{2j} G_{m(1-\varphi)}(t) \Big\|_{2}^{2}&= 
\int \Big( \big| \xi \big|^{r(m+1)+ 2mj} e^{-t|\xi|^m} 
\big( 1-\varphi \big) \Big)^2 d\xi\\
&= \int_{|\xi| \geq 1} \big| \xi \big|^{2r(m+1)+ 4mj} e^{-2t|\xi|^m} d\xi  
\leq e^{-t}\int \big| \xi \big|^{2r(m+1)+ 4mj} e^{-t|\xi|^m} d\xi  \\
&\quad\quad\leq Ce^{-t} t^{- \frac{2}{m} [r(m+1)+ 2jm]- \frac{1}{m}}\;.
\end{align*}
Hence, returning to (\ref{lem2.3.8}), we have
\begin{multline}\label{lem2.3.9}
\nor \sum_{j=0}^{N}\frac{t^j}{j!}
 K_m^j\ast \big( \partial_x M \big)^r M^{2j} G_{m(1-\varphi)}(t)
\ast v_0 \nor_2 \leq 
C\sum_{j=0}^{N}\frac{t^{j}}{j!}
e^{-\frac{t}{2}}t^{ - \frac{1}{m} [r(m+1)+ 2jm]- \frac{1}{2m} }\big\| v_0 \big\|_1 \\
\leq Ce^{-\frac{t}{2}}t^{-\frac{r(m+1)}{m} -\frac{1}{2m}}
\biggl(\sum_{j=0}^{N}\frac{t^{-j}}{j!}\biggr) \big\| v_0 \big\|_1\:.
\end{multline}
The Lemma \ref{lemab} is consequence of (\ref{lem2.3.7}) and  (\ref{lem2.3.9}).
\hfill$\BO$
%
%
\begin{lem}\label{lemad}
Let $(k,r,j)\in \N^3$. Then there exist a constant $C=C(k,r,j) >0$ 
such that
\begin{multline}\label{lemac1}
\nor \big( \partial_x M \big)^r M^{2j} G_m(t)\ast v_0 - 
\sum_{\alpha=0}^k \frac{(-1)^{\alpha}}{\alpha !}
\biggl(\int v_0 x^{\alpha}dx\biggr)
 \big( \partial_x M \big)^r M^{2j} \partial_x^{\alpha} G_m(t) \nor_2  \\
\leq Ct^{- \frac{1}{m}(r(m+1) + 2mj + k + 1) - \frac{1}{2m}}
 \big\| |x|^{k+1}v_0 \big\|_1 \;,
\end{multline}
for all $t>0$ and $v_0 \in \Le^1 \big( \R,1+|x|^{k+1} \big)$.
\end{lem}
%
%
\proof
In \cite{Duozuazua} it is proved that $v_0$ can be decomposed as
$$
v_0 = \sum_{\alpha =0}^k \frac{(-1)^{\alpha}}{\alpha !}
\biggl(\int v_0 x^{\alpha}dx\biggr) \partial_x^{\alpha}\delta +
\partial_x^{k+1} F_{k+1}\,,$$
where $F_{k+1} \in \Le^1(\R)$, such that 
$\big\| F_{k+1} \big\|_1 \leq \big\| |x|^{1+k}v_0 \big\|_1 $.
Then, taking the convolution of $v_0$ with
$\big( \partial_x M \big)^r M^{2j} G_m(x,t)$ we have
\begin{align*}
\nor \big( \partial_x M \big)^r & M^{2j} G_m(t)\ast v_0 -
\sum_{\alpha =0}^k \frac{(-1)^{\alpha}}{\alpha !}
\biggl(\int v_0 x^{\alpha}dx\biggr)
\big( \partial_x M \big)^r M^{2j} \partial_x^{\alpha}G_m(t) \nor_2 \\
&\leq
\big\| F_{k+1}\ast \big( \partial_x M \big)^r M^{2j} 
\partial_x^{k+1}G_m(t)\big\|_2 
\leq
\big\| F_{k+1} \big\|_1 
\big\| \big( \partial_x M \big)^r M^{2j} \partial_x^{k+1}G_m(t) \big\|_2\\
&\qquad\leq C\big\|\; |x|^{1+k}v_0 \big\|_1 
t^{- \frac{1}{m}(r(m+1) + 2mj + k + 1)-\frac{1}{2m}}.\hspace{2cm}{\BO}
\end{align*}
%
{{\bf Proof of Theorem \ref{theoa}}\;:\\
%
>From Lemma \ref{lemab} it follows that
\begin{equation}\label{theo2.1.1}
\nor S(t)v_0 - 
\sum_{r=0}^{N}\frac{t^r}{r!} 
K_m^r\ast \big( \partial_x M \big)^r
S_\varphi(t)\ast v_0 \nor_2  \leq 
Ct^{-\frac{N+1}{m} - \frac{1}{2m}} \big\| v_0 \big\|_1 
+ Ce^{-\frac{t}{2}} \big\| v_0 \big\|_2 \;.
\end{equation}
Moreover from Lemma \ref{lemac} we have 
\begin{multline}\label{theo2.1.2}
\nor \sum_{r=0}^{N}\frac{t^r}{r!}K_m^r\ast \big( \partial_x M \big)^r
 S_{\varphi}\ast v_0 
-\sum_{r=0}^{N}\frac{t^r}{r!}
\sum_{j=0}^{\maxi \frac{N}{2} \maxd} \frac{t^j}{j!}
K_m^{r+j} \ast \big( \partial_x M \big)^r M^{2j} G_m\ast v_0 \nor_2   \\  
\leq Ct^{-(\maxi \frac{N}{2} \maxd +1)-\frac{1}{2m}}
\biggl( \sum_{r=0}^{N} \frac{t^{-\frac{r}{m}}}{r!}\biggr)\big\| v_0 \big\|_1 
\quad+ Ce^{-\frac{t}{2}} t^{-\frac{1}{2m}}
\biggl( \sum_{r=0}^{N}\frac{t^{- \frac{r}{m}}}{r!}\biggr)
\biggl( \sum_{j=0}^{\maxi \frac{N}{2} \maxd} \frac{t^{-j}}{j!}\biggr)
\big\| v_0 \big\|_1\:. 
\end{multline}
Then from (\ref{theo2.1.1}) and (\ref{theo2.1.2}) it follows that
\begin{align}\label{theo2.1.3}
\nor S(t) v_0 \;- 
\sum_{r=0}^{N} &\frac{t^r}{r!}
\sum_{j=0}^{\maxi \frac{N}{2} \maxd} \frac{t^j}{j!}
K_m^{r+j}\ast \big( \partial_x M \big)^r
M^{2j} G_m(t)\ast v_0 \nor_2  \\ 
\leq C&t^{-\frac{N+1}{m} - \frac{1}{2m}}  \big\| v_0 \big\|_1 
+\;Ct^{-(\maxi \frac{N}{2} \maxd +1)-\frac{1}{2m}}
\biggl( \sum_{r=0}^{N} \frac{t^{-\frac{r}{m}}}{r!}\biggr)\big\| v_0 \big\|_1 \nonumber \\
+ &Ce^{-\frac{t}{2}} \big\| v_0 \big\|_2 +
Ce^{-\frac{t}{2}} t^{-\frac{1}{2m}}
\biggl( \sum_{r=0}^{N}\frac{t^{-\frac{r}{m}}}{r!}\biggr)
\biggl( \sum_{j=0}^{\maxi \frac{N}{2} \maxd} \frac{t^{-j}}{j!}\biggr)
\big\| v_0 \big\|_1 \:.\nonumber
\end{align}
On the other hand from the Lemma \ref{lemad}, we obtain
\begin{align}\label{theo2.1.4}
\nor \sum_{r=0}^{N}\frac{t^r}{r!}\sum_{j=0}^{\maxi \frac{N}{2} \maxd}
\frac{t^{j}}{j!}
K_m^{r+j}& \ast \big( \partial_x M \big)^r M^{2j} G_m(t)\ast v_0  \\
-& \sum_{r=0}^{N}\frac{t^r}{r!}
\sum_{j=0}^{\maxi \frac{N}{2} \maxd}\frac{t^{j}}{j!}
\sum_{\alpha=0}^{k(r,j)}
{\cal{M}}_{\alpha}(v_0)
K_m^{r+j}\ast \big( \partial_x M \big)^r M^{2j} 
\partial_x^{\alpha} G_m(t) \nor_2  \nonumber \\
\leq \sum_{r=0}^{N}\frac{t^r}{r!}
\sum_{j=0}^{\maxi \frac{N}{2} \maxd}\frac{t^{j}}{j!}
\nor & K_m^{r+j}\ast \big( \partial_x M \big)^r  M^{2j} G_m(t)\ast v_0  \nonumber \\
& \hspace{2cm}- \sum_{\alpha=0}^{k(r,j)}
{\cal{M}}_{\alpha}(v_0)
K_m^{r+j} \ast \big( \partial_x M \big)^r M^{2j} 
\partial_x^{\alpha} G_m(t) \nor_2 \nonumber \\
\leq C \sum_{r=0}^{N}\frac{t^r}{r!}& \sum_{j=0}^{\maxi \frac{N}{2}\maxd}
\frac{t^{j}}{j!}
\:t^{-\frac{1}{m}[k(r,j)+1+r(m+1)+2mj] - \frac{1}{2m}}\big\| |x|^{k(r,j)+1}v_0 \big\|_1
\;.\nonumber
\end{align}
In (\ref{theo2.1.4}), $k(r,j)$ denotes a natural number and we use that
$\big\| K_m^r \big\|_1 = 1$.\\
We choose $k(r,j) = N-r-mj  \geq 0$ in (\ref{theo2.1.4}). 
Hence
\begin{align}\label{theo2.1.5}
\nor & \sum_{r=0}^{N} \frac{t^r}{r!}
\sum_{j=0}^{\maxi \frac{N}{2} \maxd} \frac{t^{j}}{j!}
K_m^{r+j}\ast (\partial_x M)^r M^{r+j} G_m(t)\ast v_0 -  \\
&\quad\quad\quad - \sum_{r=0}^{N}\frac{t^r}{r!}\sum_{j=0}^{\maxi \frac{N}{2} \maxd}
\frac{t^{j}}{j!}
\sum_{\scriptstyle 0\leq \alpha \leq N-r-mj}^{*}
{\cal{M}}_{\alpha}(v_0)
K_m^j\ast \big( \partial_x M \big)^r M^{2j} 
\partial_x^{\alpha} G_m(t) \nor_2  \nonumber \\
&\quad\quad\quad\quad\quad \leq Ct^{-\frac{N+1}{m}-\frac{1}{2m}}
\sum_{r=0}^{N}
\sum_{\scriptstyle 0\leq j \leq \maxi \frac{N}{2} \maxd 
\atop\scriptstyle 0\leq j \leq \maxi \frac{N-r}{m} \maxd}
\big\| |x|^{(N+1)-r-jm}v_0 \big\|_1 \:.\nonumber
\end{align}
Finally, Theorem \ref{theoa} is consequence of
(\ref{theo2.1.3}) and (\ref{theo2.1.5}). Indeed,
\begin{align*}
\nor S(t)v_0  & -
\sum_{r=0}^{N}\frac{t^r}{r!}\sum_{j=0}^{\maxi \frac{N}{2} \maxd}
\frac{t^{j}}{j!}
\sum_{\scriptstyle 0 \leq \alpha\leq N-r-mj }^{*}
{\cal{M}}_{\alpha}(v_0) 
K_m^{r+j}\ast \big( \partial_x M \big)^r M^{2j} \partial_x^{\alpha}G_m(t) \nor_2  \\
& \leq Ct^{-\frac{N+1}{m} - \frac{1}{2m}}  \big\| v_0 \big\|_1 
+\;Ct^{-(\maxi \frac{N}{2} \maxd +1)-\frac{1}{2m}}
\biggl( \sum_{r=0}^{N} \frac{t^{-\frac{r}{m}}}{r!}\biggr)\big\| v_0 \big\|_1 + \\
&\quad+ Ct^{-\frac{N+1}{m}-\frac{1}{2m}}
\sum_{r=0}^{N}
\sum_{\scriptstyle 0\leq j \leq \maxi \frac{N}{2} \maxd 
\atop\scriptstyle 0\leq j \leq \maxi \frac{N-r}{m} \maxd}
\big\| |x|^{N+1-r-jm}v_0 \big\|_1 +\\
&\quad \quad + Ce^{-\frac{t}{2}} \big\| v_0 \big\|_2 +
Ce^{-\frac{t}{2}} t^{-\frac{1}{2m}}
\biggl( \sum_{r=0}^{N}\frac{t^{-\frac{r}{m}}}{r!}\biggr)
\biggl( \sum_{j=0}^{\maxi \frac{N}{2} \maxd} \frac{t^{-j}}{j!}\biggr)
\big\| v_0 \big\|_1 \:.\hspace{2cm} {\BO}
\end{align*}
%
{\bf Proof of Theorem \ref{theob}\;:}\\
In this case, are continued the same steps of the proof of previous theorem,
to consider $N=1$ and choose $k(r)=1-r$ in (\ref{theo2.1.4}).\hfill $\BO$
%
\setcounter{equation}{0}
\section{\hspace{-5mm}.\hspace{5mm}Complete asymptotic expansion.}
%
The function entering in the $\Le^2$-estimate of Theorem \ref{theoa}
can be written as
\begin{multline}\label{sim0}
S(t)v_0 
\;-\;\sum_{\alpha =0}^{N}
{\cal{M}}_{\alpha}(v_0) \partial_x^{\alpha} G_m(t) \\
-\;\sum_{r=0}^{N}\frac{t^r}{r!}\sum_{j=0}^{\maxi \frac{N}{2} \maxd}
\frac{t^{j}}{j!}
\sum_{\scriptstyle 0\leq \alpha \leq N-r-mj \atop\scriptstyle (r,j)\not= (0,0)}^{*} 
{\cal{M}}_{\alpha}(v_0)\;
K_m^{r+j} \ast \big(\partial_x M \big)^r M^{2j} \partial_x^{\alpha} G_m(t)\:.
\end{multline}
In this section we simplify the terms of (\ref{sim0}) that have the form
\begin{equation}\label{sim1}
K_m^{r+j}\ast \big( \partial_x M \big)^r M^{2j} \partial_x^{\alpha} G_m(t),
\end{equation}
and we replace them, for the first term of their asymptotic expansion, i.e. by
$$
(\partial_x M)^r M^{2j} \partial_x^{\alpha} G_m(t).
$$
Now, of the moments formula
$$ 
\int_{\R} x^{\alpha} f(x)dx = 
i^{\alpha}\left( \partial_{\xi}^{\alpha}\fo{f}\right)(0)\;,\: \alpha>0.
$$
Then, observing that 
\begin{equation}\label{momento1}
\hspace{3cm}\left\{
\begin{array}{ll}
\int_{\R} K_m^j(x)dx = 1,  \\
\int_{\R} x K_m^j(x)dx = 0\;,\quad j \geq 1, m > 1.
\end{array}
\right.
\end{equation}
On the other hand, as $K_m^j \in \Le^1(\R)$, the solution of 
\begin{equation}\label{cal}
\hspace{3cm}\left\{
\begin{array}{ll}
u_t - Mu = 0 \\
u(x,0) = K_m^j(x)\;,\;\; j \geq 1,
\end{array}
\right.
\end{equation}
is given by
$$
u(x,t) = \big( K_m^j\ast G_m(t) \big)(x)\;.
$$
Keeping (\ref{momento1}) in mind, we have
$K_m^j \in \Le^1 \big( \R,(1+|x| \big)$, since
$$
\int\big| x \big| \big| K_m^j(x) \big|dx < \infty .
$$
Consequence of Theorem \ref{calorg} and
as ${\cal{M}}_0 (K_m^j)= \int K_m^j =1$,
we have the following Corollary:
%
\begin{cor}\label{sim3}
There exist a constant $C=C(m)>0$ such that
$$
\big\| G_m(t)\ast K_m^j - G_m(t) \big\|_2
\leq C\;t^{-\frac{1}{m} - \frac{1}{2m}}
\big\| |x| K_m \big\|_1 \:.
$$
\end{cor}
\proof
By Lemma (\ref{lemaa}) with $N=0$ and $\int K_m(x)dx=1$
we have
\begin{equation}\label{prev0}
\big\| G_m(t)\ast K_m - G_m(t) \big\|_2 \leq C\;t^{-\frac{1}{m} - \frac{1}{2m}}
\big\| |x| K_m \big\|_1\;.
\end{equation}  
Then, from (\ref{prev0}) we have
\begin{align}\label{prev1}
\big\| G_m(t)\ast K_m^2 - G_m(t)\ast K_m \big\|_2 &\leq
\big\| K_m \big\|_1 \big\| G_m(t)\ast K_m - G_m(t) \big\|_2  \\
& \leq 
C\;t^{-\frac{1}{m} - \frac{1}{2m}}
\big\| |x| K_m \big\|_1\;.\nonumber
\end{align}
Now, from (\ref{prev0}) and (\ref{prev1}), it follows that
\begin{multline*}
\big\| G_m(t)\ast K_m^2 - G_m(t) \big\|_2 \leq 
\big\| G_m(t)\ast K_m^2 - G_m(t)\ast K_m \big\|_2  \\
+\; \big\| G_m(t)\ast K_m - G_m(t) \big\|_2
\leq C\;t^{-\frac{1}{m} - \frac{1}{2m}} \big\| |x| K_m \big\|_1\;.
\end{multline*}
The conclusion of the proof the Corollary \ref{sim3} it follows by 
induction.
\hfill $\BO$\\
If instead of $G_m(t)$ we consider 
$(\partial_x M)^r M^{2j} \partial_x^{\alpha}G_m(t)$
in the previous Corollary we obtain the following result:
%
\begin{cor}\label{sim4}
Let $r\in \N $. Then there exist a constant $C=C(r,j)>0$ such that
\begin{multline}
\Big\| (\partial_x M)^r M^{2j} \partial_x^{\alpha}G_m (t)\ast K_m^{r+j}\; -\;
(\partial_x M)^r M^{2j} \partial_x^{\alpha}G_m(t) \Big\|_2 \leq \\
\leq Ct^{-\frac{1+(m+1)r+2jm+|\alpha|}{m} - \frac{1}{2m}}
\big\| |x| K_m \big\|_1 \;.
\end{multline}
\end{cor}
%
And, as an immediate consequence of Corollary {\ref{sim4} we have:
%
\begin{cor}\label{sim5} 
Let $r\in \N $.  
Then there exist a constant $C=C(r,j)>0$ such that
\begin{align*}
\nor
\sum_{r=0}^{N}\frac{t^r}{r!}
\sum_{j=0}^{\maxi \frac{N}{2} \maxd}\frac{t^j}{j!}
&\sum_{\scriptstyle 0\leq \alpha \leq N-r-mj \atop\scriptstyle (r,j)\not= (0,0)}^{*}
\hspace{-0.3cm} {\cal{M}}_{\alpha}(v_0)
\biggl(K_m^{r+j} \ast \big( \partial_x M \big)^r M^{2j} \partial_x^{\alpha}G_m(t) 
- \big( \partial_x M \big)^r M^{2j} \partial_x^{\alpha}G_m(t) \biggr)
\nor_2  \\
&\leq 
C\;\sum_{r=0}^{N}\frac{t^r}{r!}
\sum_{j=0}^{\maxi \frac{N}{2} \maxd}\frac{t^j}{j!}
\sum_{0 \leq \alpha\leq N-r-mj}^{*} |{\cal{M}}_{\alpha}(v_0)|
\;t^{-\frac{1+(m+1)r+2jm+\alpha}{m} - \frac{1}{2m}}
\big\| |x| K_m \big\|_1  \\
&\quad\quad\leq C
\max_{0\leq \alpha \leq N-1}\biggl\{|{\cal{M}}_{\alpha}(v_0)|\biggr\}
\;t^{-\frac{1+N}{m} - \frac{1}{2m}}
\big\| |x|K_m \big\|_1\:.
\end{align*}
\end{cor}
%
Then the proof of Theorem \ref{theo1} follows directly from Corollary \ref{sim5}
and equation (\ref{sim0}), indeed:
\begin{align*}
\nor S(t)v_0&\;-\; 
\hspace{-0.2cm}\sum_{0\leq\alpha\leq N} {\cal{M}}_{\alpha}(v_0) \partial_x^{\alpha}G_m(t) 
\;-\; \sum_{r=0}^{N}\frac{t^r}{r!}
\sum_{j=0}^{\maxi \frac{N}{2} \maxd}\frac{t^{j}}{j!}
\sum_{\scriptstyle 0\leq \alpha \leq N-r-mj \atop\scriptstyle (r,j)\not= (0,0)}^{*}
\hspace{-0.3cm}{\cal{M}}_{\alpha}(v_0)
\big( \partial_x M \big)^r M^{2j} \partial_x^{\alpha}G_m(t) 
\nor_2\\
&\leq C \Bigg[
t^{- (\maxi \frac{N}{2} \maxd + 1) -\frac{1}{2m}} \big\| v_0 \big\|_1  
\biggl( \sum_{r=0}^{N} \frac{t^{-\frac{r}{m}}}{r!}\biggr)
+t^{-\frac{N+1}{m}-\frac{1}{2m}}\big\| v_0 \big\|_1
+t^{-\frac{N+1}{m}-\frac{1}{2m}}\sum_{r=1}^{N+1}
\big\| |x|^r v_0 \big\|_1   \\ 
& \hspace{4cm}+ \,e^{-\frac{t}{2}} \big\| v_0 \big\|_2\, +\,
e^{-\frac{t}{2}} t^{-\frac{1}{2m}}
\biggl( \sum_{r=0}^{N}\frac{t^{-\frac{r}{m}}}{r!}\biggr)
\biggl( \sum_{j=0}^{\maxi \frac{N}{2} \maxd}
\frac{t^{-j}}{j!}\biggr) \big\| v_0 \big\|_1  \\
&\hspace{6cm}+\;\max_{0\leq \alpha \leq N-1} \biggl\{|{\cal{M}}_{\alpha}(v_0)|\biggr\}
\;t^{-\frac{1+N}{m} - \frac{1}{2m}} \big\| |x| K_m \big\|_1 
\Bigg]
\hspace{0.5cm}{\BO}
\end{align*}
%
\begin{remark}
When $m=2$ the Theorem \ref{theo1} it agree with results in \cite{razu}\\
Remembering that
${\cal{M}}_{\alpha}(v_0)\,=\, 
\frac{(-1)^{\alpha}}{\alpha !}\int x^{\alpha}v_0(x)dx$. In the Theorem \ref{theo1}
we have the following particular cases when $m \in Z^{+}$, $m>1$:
\begin{enumerate}
\item 
If $N=0$, 
the firs term in the asymptotic expansion is  
$r_1 (x,t) = {\cal{M}}_0 (v_0)\;G_m(t)$. Then, for
$v_0\in \Le^2(\R) \cap \Le^1 (\R,1+|x|)$, we have 
$$
t^{\frac{1}{2m}} \big\| S(t)v_0 -{\cal{M}}_0(v_0) G_m(t) \big\|_2 \leq 
Ct^{-\frac{1}{m}}\,, 
\quad \mbox{for}\quad t \geq 1\,.
$$
%
\item 
If $N=1$, 
the firs term in the asymptotic expansion is 
$r_1 (x,t) =  {\cal{M}}_0 (v_0) G_m(t)$, 
and the secon term is  
$r_2 (x,t) = {\cal{M}}_1 (v_0) \partial_x G_m(t) +
t {\cal{M}}_0 (v_0) (\partial_{x} M) G_m(t)$.\\
Then for
$v_0 \in \Le^2 (\R)\cap \Le^1(\R,1+|x|^2)$ we have 
$$
t^{\frac{1}{2m}+\frac{1}{m}} \nor S(t)v_0 - \sum_{j=1}^{2} r_j (t) \nor_2 
\leq Ct^{-\frac{1}{m}}\,,\quad \text{for} \quad t\geq 1\,.$$
%
\item 
If $N=2$, $m=2$ then, 
$Mu = -\partial_{xx}u$ and 
$G_2(x,t)= G(x,t)=(4\pi t)^{-\frac{1}{2}}e^{-\frac{x^2}{4t}}$ 
it is the heat kernel. The firs term in the asymptotic expansion is  \\
$r_1 (x,t) ={\cal{M}}_0 (v_0)G(t)$, the secon term is \\
$r_2 (x,t) = 
{\cal{M}}_1 (v_0) \partial_x G(t) - t{\cal{M}}_0 (v_0)\partial_{xxx}G(t)$.
And the third one is\\
$r_3 (x,t) = {\cal{M}}_2 (v_0)\partial_{xx} G +
t\Big( {\cal{M}}_0 (v_0) - {\cal{M}}_1 (v_0)\Big)\partial_x^4 G
+\frac{t^2}{2}{\cal{M}}_0 (v_0)\partial_x^6 G.$ \\
Then for
$v_0 \in \Le^2 (\R)\cap \Le^1(\R,1+|x|^3)$ we have 
$$
t^{\frac{1}{4}+1} \nor S(t)v_0 - \sum_{j=1}^{3} r_j (t) \nor_2 
\leq Ct^{-\frac{1}{2}}\,,\quad \text{for} \quad t\geq 1.
$$.\\
When $m\geq 3$, $N=2$, the firs term in the asymptotic expansion is  \\
$r_1 (x,t) ={\cal{M}}_0 (v_0)G_m(t)$, the secon term is \\
$r_2 (x,t) = 
{\cal{M}}_1 (v_0) \partial_x G_m(t) + t{\cal{M}}_0 (v_0)(\partial_{x}M)G_m(t)$.
And the third one is\\
$r_3 (x,t) = {\cal{M}}_2 (v_0)\partial_{xx} G_m(t) +
t{\cal{M}}_1 (v_0)(\partial_x M) \partial_x G_m(t)
+\frac{t^2}{2}{\cal{M}}_0 (v_0)(\partial_x M)^2 G_m(t).$\\
Then for
$v_0 \in \Le^2 (\R)\cap \Le^1(\R,1+|x|^3)$ we have 
$$
t^{\frac{1}{2m}+\frac{2}{m}} \nor S(t)v_0 - \sum_{j=1}^{3} r_j (t) \nor_2 
\leq Ct^{-\frac{1}{m}}\,,\quad \text{for} \quad t\geq 1,\quad m \geq 3.$$
\end{enumerate}
\end{remark}
\begin{remark}
When $m = n+\delta$ for $n \in Z^{+}$, $n>1$, $0 < \delta < 1$, as consequence 
of Theorem \ref{theo1}, we have
$$
t^{\frac{1}{m} +\frac{1}{2m}}
\Big\| S(t)v_0  -{\cal{M}}_0(v_0)G_m(t)- {\cal{M}}_1(v_0)
\partial_x G_m(t)- t{\cal{M}}_{0}(v_0) \big( \partial_x M \big)G_m(t) \Big\|_2 \leq 
C t^{-\frac{1}{m}},
$$
for $t \geq 1$ and
for all $v_0\in \Le^2(\R) \cap \Le^1(\Rn;1+|x|^2)$.
\end{remark}

%
\setcounter{equation}{0}
\section{\hspace{-5mm}. Complete asymptotic expansion of the   
generalized KdV linear equation}
%
Let us consider now the linearized equation
\begin{equation}\label{kdvg1}
\left\{ \begin{array}{rcl}
u_t + Mu - Mu_x &=& 0,
\quad x \in \R, t > 0\\
u(x,0) &=&u_0(x)\,.
\end{array}\right. 
\end{equation}
Where $\fo{Mu}(\xi)= |\xi|^m \fo{u}(\xi)$, $m\in \Z^{+}$.
Its solution takes the form
$$u(x,t) = T(\cdot,t) \ast u_0(x)\,,$$
with
$$T(x,t) = \frac{1}{2\pi} \int e^{t\psi(\xi) + ix\cdot\xi}
d\xi.$$
The phase function is now 
$$\psi (\xi) = -|\xi|^m + i\xi |\xi|^m$$ which is
a bit simpler than the phase function $\widetilde{\Phi}$ of the (\ref{genelin}) linear
equation.
Proceeding as in the proof of Theorem \ref{theoa},
Section 3 we obtain the following result:
\begin{theo}\label{theokdvg}
For any $N \in \N$, there exists a constant $C=C(N)>0$ such that
\begin{multline*}
\nor T(t)\ast u_0 - \sum_{j=0}^{N} \frac{t^j}{j!} \big( \partial_x M \big)^j
\sum_{\alpha = 0}^{N-j} {\cal{M}}_{\alpha} (u_0) \partial_x^{\alpha} G_m(t)\nor_2 \\
\leq C\,t^{-\frac{1}{2m}-\frac{N+1}{m}} \big\| u_0 \big\|_1 + 
C\, t^{-\frac{1}{2m}-\frac{N+1}{m}} 
\sum_{k=1}^{N+1} \big\| |x|^k u_0 \big\|_1\,,
\end{multline*}
for all $t>0$ and $u_0 \in \Le^1 (\R,1+|x|^{N+1})$.
\end{theo}
\begin{remark}
For instance, if $N=2$, in Theorem \ref{theokdvg} we have that
the first term is $r_1 (x,t) = MG_m(t)$.
The second term is
$r_2 (x,t) = {\cal{M}}_1 (u_0) G_m(t) +
t {\cal{M}}_0 (u_0) \big( \partial_x M\big) G_m(t)$,
and the third one
$$
r_3 (x,t) ={\cal{M}}_2 (u_0)\partial_{xx}G_m(t)+
t{\cal{M}}_1 (u_0) \big( \partial_x M \big) \partial_x G_m(t)
+ \frac{t^2}{2} {\cal{M}}_0 (u_0)  \big( \partial_x M \big)^2 G_m(t).
$$
Then for
$u_0 \in \Le^1 (\R,1+|x|^3)$ we have
$$t^{\frac{1}{4}+1} \nor T(t)u_0 - \sum_{j=1}^{3} r_j (t) \nor_2 \leq 
Ct^{-\frac{1}{m}}\,,\quad \mbox{for} \quad t\geq 1\,.$$
\end{remark}
%
\begin{remark}
In accordance with Theorems \ref{theo1} and \ref{theokdvg}, the succesive
terms appearing in the asymptotic expansion of the solutions of (\ref{genelin})
and (\ref{kdvg1}) have the form 

\begin{equation}\label{ter1}
\sum_{r=0}^{N}\frac{t^r}{r!}
\sum_{j=0}^{\maxi \frac{N}{2} \maxd}\frac{t^{j}}{j!}
\big( \partial_x M \big)^r M^{2j}
\sum_{\alpha=0}^{N-r-2j}
{\cal{M}}_{\alpha}(v_0) \partial_x^{\alpha}G_m(t)
\end{equation}
and
\begin{equation}\label{ter2}
\sum_{k=0}^{N} \frac{t^k}{k!}\;\big( \partial_x M \big)^k
\sum_{\alpha =0}^{N-k} {\cal{M}}_{\alpha} (u_0) \partial_x^{\alpha} G_m(t)\;,
\end{equation}
respectively. We see that, the term due to dispersive effects of 
$Mu_x$ in (\ref{kdvg1}) and 
$Mu_x$ and $Mv_t$ in (\ref{genelin}),
appear in the asymptotic expansions starting at the second term.
\end{remark}
%
%
\section{Global Solution}
%
In this section we study the global solution of the initial value problem for 
the following model equation (\ref{gene1}) and suppose that $m>2$.
There is a well known principle which has frequently been used to prove 
existence of global solution of non-linear equations. 
Indeed, we may rewrite our non-linear diferential partial equation as 
non-linear equation integral, obtained from the formula of variation of the 
parameters (or formula of Duhamel)
\begin{gather}\label{intequa}
v(x,t) = 
S(t)v_0(x) - \int_0^t S(t-\tau) K_m \ast \big( v^q \big)_x (\tau)d\tau.
\end{gather}
Where
\begin{equation}\label{solinear2}
S(t)v_0(x) = \frac{1}{(2\pi)} 
\int \mbox{exp} \big( -t\Phi(\xi) + ix\cdot\xi \big) \fo{v_0} (\xi) d\xi\,,
\end{equation}
with the phase function
$$
\Phi(\xi)= \frac{|\xi|^m - i\xi|\xi|^m}{1 + |\xi|^m}.
$$
Recall that $K_m$ is the function defined by
$\fo{K}(\xi)=\frac{1}{1+|\xi|^m}$.\\
The equation (\ref{intequa}) is equivalent to the differential form (\ref{gene1}), but is much easier to 
handle when it comes to proving questions of existence an uniqueness.

To find a solution $v$ to (\ref{intequa}), we shall use an iterative method. 
We first approximate $v$ by the linear solution
\begin{gather}\label{linsol}
v_0(x,t) = S(t)v_0(x).
\end{gather}
Then we make a better approximation
$$
v_1(x,t) = 
S(t)v_0(x) - \int_0^t S(t-\tau) K_m \ast \big( v_0^q \big)_x (\tau)d\tau.
$$
More generally, we define the non-linear map $v \mapsto Pv$ by
$$
Pv(x,t) = 
S(t)v_0(x) - \int_0^t S(t-\tau) K_m \ast \big( v^q \big)_x (\tau)d\tau
$$
and define $v_{k+1}=P(v_k)$ for all $k=0,1,...$. We hope to show that this
sequence of approximations converges to a limit $v$, so that $v=Pv$.
This would give us a solution to (\ref{intequa}).

In short, we want to find a fixed point of $P$, show that it is unique.
We can accomplish all this in one stroke from the Banach´s fixed point Theorem,
as soon as we show that $P$ is a contraction on some complete metric space 
$X$ which contains $v_0$. 
Well, we have to pick the right complete metric space to get the contraction 
working.
\begin{theo}\label{global}
Suppose that $v_0 \in W^{2,1}(\R)$; then there exist a positive constant 
$\delta_1$ such that when $\|v_0\|_{W^{2,1}} < \delta_1$, then 
the equation (\ref{intequa}), have a unique global solution $v(x,t)$ 
satisfying
$$
v(x,t)\in C(0,\infty; \Le^{\infty}(\R) \cap H^2(\R) ).
$$
Moreover, the asymptotic decay rates of the solutions $v(x,t)$, (\ref{optimos}),
holds.
\end{theo}
Initial data $v_0$ for which $\|v_0\|_{W^{2,1}}$ is sufficiently small give 
rise to global solution of the nonlinear equation.
\subsection{Properties of the linear solutions $S(t)v_0$}
%
\begin{lem}\label{svde}
Let $v_0 \in W^{2,1}(\R)$. Then there exist a positive constant $C$, independent
of $t$ such that
\begin{align}
\big\|& S(t)v_0 \big\|_2 
\leq C\|v_0\|_{W^{1,1}} (1+t)^{-\frac{1}{2m}} \label{sv}\\
\big\|& S_x(t)v_0 \big\|_2 \leq 
C\|v_0\|_{W^{2,1}}(1+t)^{-\frac{1}{2m}- \frac{1}{m}}  \label{sxv}\\
\big\|& S(t)v_0\big\|_{\infty} \leq 
C\|v_0\|_{W^{2,1}} (1+t)^{-\frac{1}{m}}. \label{svi}
\end{align}
\end{lem}
%
\proof
By Plancherel's formula, we have
\begin{gather}\label{sv1}
\big\| S(t)v_0 \big\|_2^2 =
\int \frac{ \mbox{e}^{-\frac{2t|\xi|^m}{1+|\xi|^m}} }{1+|\xi|^2}
(1+|\xi|^2) |\fo{v_0}(\xi)|^2 d\xi 
\leq \bigg[\sup_{\xi \in \R} \Big\{ (1+|\xi|) |\fo{v_0}(\xi)|\Big\} \bigg]^2 
\int_{\R} \frac{ \mbox{e}^{-\frac{2t|\xi|^m}{1+|\xi|^m}} }{1+|\xi|^2} .
\end{gather}
Now
\begin{gather}\label{in}
\int_{\R} \frac{ \mbox{e}^{-\frac{2t|\xi|^m}{1+|\xi|^m}} }{1+|\xi|^2}
\leq 2\bigg( \int_{0}^{1} + \int_{1}^{\infty} \bigg).
\end{gather}
If $\xi \in [0,1]$ then $\frac{1}{2} \leq \frac{1}{1+|\xi|^m} \leq 1$, hence 
$\mbox{e}^{-\frac{2t|\xi|^m}{1+|\xi|^m}} \leq \mbox{e}^{-t|\xi|^m}$, we have
\begin{align}\label{sv2}
\int_0^1 \frac{ \mbox{e}^{\frac{-2t|\xi|^m}{1+|\xi|^m}} }{1+|\xi|^2} d\xi 
 \leq& C \int_0^1 \mbox{e}^{-t|\xi|^m} d\xi 
\leq C \int_0^1 \mbox{e}^{-(1+t)|\xi|^m}\mbox{e}^{|\xi|^m} d\xi \\
&\leq  C \int_0^1 \mbox{e}^{-(1+t)|\xi|^m}d\xi
\leq C (1+t)^{-\frac{1}{m}}. \nonumber
\end{align}
For the second term on the right-hand side of (\ref{in}), if 
$\xi \in [1,\infty)$ then $\frac{1}{2} \leq \frac{|\xi|^m}{1+|\xi|^m} \leq 1$, hence
\begin{gather}\label{sv3}
\int_1^{\infty} \frac{ \mbox{e}^{ \frac{-2t|\xi|^m}{1+|\xi|^m} } }{1+|\xi|^2} d\xi 
\leq \int_1^{\infty} \frac{ \mbox{e}^{-t } }{1+|\xi|^2} d\xi
\leq C \mbox{e}^{-t}.
\end{gather}
On the other hand 
\begin{align}\label{sup1}
\sup_{\xi \in \R} \Big\{ (1+|\xi|) |\fo{v_0}(\xi)| \Big\} 
&\leq \sup_{\xi \in \R} \Bigg| \int e^{-ix\cdot\xi} v_{0}(x) dx \Bigg|
+  \sup_{\xi \in \R} \Bigg| \int e^{-ix\cdot\xi} v_{0x}(x) \Bigg| dx \\
&\quad \leq \int |v_{0}(x)|dx +  \int |v_{0x}(x)|dx \leq \|v_{0} \|_{W^{1,1}}.\nonumber
\end{align}\\
Thus, plugging (\ref{sv2}), (\ref{sv3}) and (\ref{sup1}) into (\ref{sv1}) implies 
({\ref{sv}}).

{\it{Proof of}} (\ref{sxv}), we have similarly
\begin{multline}\label{sxv1}
\big\| S_x(t)v_0 \big\|_2^2 =
\int  
\frac{(1+|\xi|^2)^2 |\xi|^2 \mbox{e}^{-\frac{2t|\xi|^m}{1+|\xi|^m}} }{(1+|\xi|^2)^2} 
|\fo{v_0}(\xi)|^2 d\xi \\
\leq \bigg[ \sup_{\xi \in \R} \Big\{ (1+|\xi|^2) |\fo{v_0}(\xi)| \Big\} \bigg]^2
\int_{\R} \frac{ |\xi|^2 \mbox{e}^{-\frac{2t|\xi|^m}{1+|\xi|^m}} }{(1+|\xi|^2)^2}.
\end{multline}
Now, 
$$
\int_{\R} \frac{ |\xi|^2 \mbox{e}^{-\frac{2t|\xi|^m}{1+|\xi|^m}} }{(1+|\xi|^2)^2}
\leq 2\bigg( \int_{0}^{1} + \int_{1}^{\infty} \bigg).
$$
Then, $\xi \in [0,1]$ implies that $\frac{1}{2} \leq \frac{1}{1+|\xi|^m} \leq 1$, we get
\begin{equation}\label{sxv2}
\int_0^1\!\! |\xi|^2 \mbox{e}^{\frac{-2t|\xi|^m}{1+|\xi|^m}} d\xi 
 \leq C\! \int_0^1\!\! |\xi|^2 \mbox{e}^{-(1+t)|\xi|^m}  \mbox{e}^{|\xi|^m} d\xi 
\leq C\! \int_0^1\!\! |\xi|^2 \mbox{e}^{-(1+t)|\xi|^m} d\xi 
\leq C (1+t)^{-\frac{2}{m} - \frac{1}{m}}. 
\end{equation}
Using the facts that $\frac{1}{2} \leq \frac{|\xi|^m}{1+|\xi|^m} < 1$ if
$\xi \in [1,\infty)$, we have
\begin{gather}\label{sxv3}
\int_1^{\infty} \frac{|\xi|^2 \mbox{e}^{-\frac{2t|\xi|^m}{1+|\xi|^m}} }{(1+|\xi|^2)^2}d\xi
\leq \int_1^{\infty} \frac{ \mbox{e}^{-t}}{ 1+|\xi|^2 } d\xi
\leq C \mbox{e}^{-t}.
\end{gather}
Also,
\begin{align}\label{sup2}
\sup_{\xi \in \R} \Big\{ (1+|\xi|^2 ) |\fo{v_0}(\xi)| \Big\} 
&\leq  \sup_{\xi \in \R} \Bigg| \int e^{-ix\cdot\xi} v_{0}(x) dx \Bigg|
+  \sup_{\xi \in \R} \Bigg| \int e^{-ix\cdot\xi} v_{0xx}(x) \Bigg| dx \\
&\quad \leq \int |v_{0}(x)|dx +  \int |v_{0xx}(x)|dx \leq \|v_{0} \|_{W^{2,1}}.\nonumber
\end{align}
Therefore (\ref{sxv2}), (\ref{sxv3}) and (\ref{sup2}) give us the desired estimate (\ref{sxv}). 

{\it{Proof of}} (\ref{svi}), using the classical inequality
$$
\big\| S(t)v_0\big\|_{\infty}\leq \big\|S(t)v_0 \big\|_2^{1/2} \big\|S_x(t)v_0 \big\|_2^{1/2},
$$
(\ref{sv}) and (\ref{sxv}), we get (\ref{svi}).\hfill$\BO$\\
Now, recall that $K_m$ is the function defined by
$\fo{K_m}(\xi)=1/(1+|\xi|^m)$, we have:
\begin{lem}\label{skm}
Let $m> 3/2$. If 
$$
S(t)K_m(x) = \frac{1}{(2\pi)}
\int \mbox{exp} \big( -t\widetilde{\Phi}(\xi) + ix\cdot\xi \big) \fo{K_m}(\xi) d\xi, 
$$
with the phase function
$$
\widetilde{\Phi}(\xi)= \frac{|\xi|^m - i\xi|\xi|^m}{1 + |\xi|^m}.
$$
Then, there exist a positive constant $C$, independent of $t$, such that
\begin{align}
\big\|& S(t)K_m \big\|_2 \leq C(1+t)^{-\frac{1}{2m}}; \label{sk}\\
\big\|& S_x(t)K_m \big\|_2 \leq C(1+t)^{-\frac{1}{2m}- \frac{1}{m}}; \label{sxk}\\
\big\|& S(t)K_m \big\|_{\infty} \leq C(1+t)^{-\frac{1}{m}}. \label{ski}
\end{align}
\end{lem}
\proof
By Plancherel formula we have
\begin{align*}
\big\| S(t)K_m \big\|_2^2 =& \big\| \fo{S(t)K_m} \big\|_2^2 =
\int \frac{ \mbox{e}^{-\frac{2t|\xi|^m}{1+|\xi|^m}} }{(1+|\xi|^m)^2} d\xi 
=\int_0^1 +\int_1^{\infty}\\
&\leq \int_0^1 \frac{ \mbox{e}^{-t|\xi|^m} }{(1+|\xi|^m)^2} d\xi +
\int_1^{\infty}\frac{\mbox{e}^{-t}}{(1+|\xi|^m)^2}d\xi
\leq C(1+t)^{-1/m} + Ce^{-t}.
\end{align*}
This  inequality implies (\ref{sk}). \\
Similarly, we have
\begin{multline*}
\big\| S_x(t)K_m \big\|_2^2 = \big\| \fo{S_x(t)K_m} \big\|_2^2 =
\int \frac{ |\xi^2| \mbox{e}^{-\frac{2t|\xi|^m}{1+|\xi|^m}} }{(1+|\xi|^m)^2} d\xi 
=\int_0^1 +\int_1^{\infty} \\
\leq \int_0^1 \frac{ |\xi|^2 \mbox{e}^{-t|\xi|^m} }{(1+|\xi|^m)^2} d\xi +
\int_1^{\infty} \frac{\mbox{e}^{-t}\;\;|\xi|^2}{(1+|\xi|^m)^2}d\xi
\leq C(1+t)^{-\frac{2}{m}- \frac{1}{m}} + Ce^{-t},\;\;(m>3/2).
\end{multline*}
Hence, (\ref{sxk}) it is proven.\\
Now,
\begin{gather*}
\big| S(t)K_m \big| \leq 
C\int \frac{ \mbox{e}^{-\frac{t|\xi|^m}{1+|\xi|^m}} }{1+|\xi|^m} d\xi 
\leq C\int_0^1 \mbox{e}^{-\frac{t|\xi|^m}{2}} d\xi +
C\int_1^{\infty} \frac{\mbox{e}^{-\frac{t}{2}} }{1+|\xi|^m}d\xi
\leq C(1+t)^{-1/m} + Ce^{-\frac{t}{2}}.
\end{gather*}
Then, (\ref{ski}) it is proven. \hfill$\BO$
\subsection{Proof of Theorem \ref{global}}
%
Looking at the behaviour of the first few iterates $v_0$, $v_1$, etc., by the 
Lema \ref{svde}, we decide that the correct metric space to use is 
(as in \cite{Cazweiss,Mei} for example)
$$
X_{\delta} = 
\Big\{ v \in C \big( 0,\infty; H^2 \cap \Le^{\infty} \big)\; 
\big/ M(v)<\delta \Big\}, \qquad m>2\;.
$$
Endowed with the distance
$$
M(v)= \sup_{0\leq t < \infty} \Big\{ 
(1+t)^{\frac{1}{2m}} \big\| v(t) \big\|_2 
+ (1+t)^{\frac{1}{m}} \big\| v(t) \big\|_{\infty} 
+(1+t)^{\frac{1}{2m}+\frac{1}{m}} \big\| v_x(t) \big\|_2  \Big\},
$$
$X_{\delta}$ is a complete metric space.
We will show that the mapping defined formally by
$$
P(v)= S(t)v_0(x) - \int_0^t S(t-\tau) K_m \ast \big( v^q \big)_x (\tau)d\tau,
$$
is a strict contraction on $X_{\delta}$, that is, 
we will prove that there exist the positive constant $\delta_1$, such that the operator
$P$ maps $X_{\delta_1}$ into itself and has a unique fixed point in 
$X_{\delta_1}$. Thus, such a fixed point $v(x,t)$ is the unique solution 
of equation (\ref{intequa}) in $X_{\delta_1}$, globally in time. To prove 
these, the following Lemmas are avaliable.
%
\begin{lem}\label{gronw}
Suppose that $a>0$ and $b>0$ and $\max \{a,b\} > 1$, then
\begin{gather}\label{gronw1}
\int_0^t (1+t-\tau)^{-a} (1+\tau)^{-b} d\tau \leq (1+t)^{-\min\{ a,b\}}.
\end{gather}
\end{lem}
%
The Lemma \ref{gronw} can be found in \cite{Segal}, see also \cite{Matsumura}.
%
\begin{lem}\label{holder}
Let $q>1$. Then there exists a constant $C>0$, independent of $t$, such that
\begin{gather}\label{holder1}
\big\| |u(t)|^{q-1}u(t) -|v(t)|^{q-1}v(t) \big\|_r \leq C \big\|u(t)-v(t)\big\|_{p_1}
\biggl( \big\|u(t)\big\|_{(q-1)p_2}^{q-1} +\big\|v(t)\big\|_{(q-1)p_2}^{q-1}\biggr),
\end{gather}
for $\frac{1}{r}=\frac{1}{p_1}\,+\, \frac{1}{p_2}\leq 1$.
\end{lem}
%
\proof
If $\alpha >0$,  we obtain
$$|a|^{\alpha}a- |b|^{\alpha}b =
(1+\alpha)(a-b)\int_{0}^{1}|\tau(a-b)+b|^{\alpha}d\tau,$$
for all $a,b \in \R$.
In our case we have that
$$|u(t)|^{\alpha}u(t)- |v(t)|^{\alpha}v(t) =
(1+\alpha)(u(t)-v(t))\int_{0}^{1}|\tau(u(t)-v(t))+v(t)|^{\alpha}d\tau\;.$$
Using the Holder inequality for
$\frac{1}{r}=\frac{1}{p_1}\,+\,\frac{1}{p_2} \leq 1$, we obtain
\begin{align*}
\|\, |u(t)|^{q-1}u(t) - |v(t)|^{q-1} v(t)\|_r &\leq C\|u(t)-v(t)\|_{p_1}
\int_{0}^{1}\|(\tau(u-v)+v)^{q-1}\|_{p_2}d\tau \\
&\leq  C\|u(t)-v(t)\|_{p_1}
\int_{0}^{1}\|\tau(u-v)+v\|_{(q-1)p_2}^{q-1}d\tau \\
&\leq C\|u(t)-v(t)\|_{p_1}
\biggl(\|u(t)\|_{(q-1)p_2}^{q-1} +\|v(t)\|_{(q-1)p_2}^{q-1}\biggr)\;.\quad{\BO}
\end{align*}
%
\begin{lem}\label{svX}
Let $q>m$ and $m>2$. If $v\in X_{\delta}$, we have 
\begin{align}
\int_0^t \Big\|& S(t-\tau) K_m \ast \big( v^q \big)_x (\tau) \Big\|_2\; d\tau 
\leq C\delta^{q}(1+t)^{-\frac{1}{2m}};    \label{svX1}\\
\int_0^t \Big\|& S_x(t-\tau) K_m \ast \big( v^q \big)_x (\tau) \Big\|_2\; d\tau 
\leq C\delta^{q}(1+t)^{-\frac{1}{2m} -\frac{1}{m}};  \label{svX2}\\
\int_0^t \Big\|& S_x(t-\tau) K_m \ast \big( v^q \big)_x (\tau) \Big\|_{\infty}\; d\tau 
\leq C\delta^{q}(1+t)^{-\frac{1}{m}}. \label{svX3}
\end{align}
\end{lem}
\proof
Note that
\begin{align}\label{svX0}
\big\| v^{(q-1)} v_x (t) \big\|_1 
\leq&  \big\| v^{(q-1)}(t)\big\|_2  \big\|v_x(t) \big\|_2 
\leq \big\| v(t) \big\|_{\infty}^{q-m} \big\| v(t) \big\|_{2(m-1)}^{m-1}  
\big\|v_x (t) \big\|_2 \\
&\leq \delta^q \; 
(1+\tau)^{-\frac{q-m}{m}-\frac{m-1}{m}(1-\frac{1}{2(m-1)})-\frac{1}{2m}-\frac{1}{m}}
\leq \delta^q (1+\tau)^{-\frac{q}{m}}.\nonumber
\end{align}
{\it{Proof of}} (\ref{svX1}),  by (\ref{sk}) and (\ref{svX0}) we have
\begin{align*}
\int_0^t \Big\| S(t-\tau) K_m \ast \big( v^q \big)_x (\tau) \Big\|_2\; d\tau 
\leq&  \int_0^t \big\| S(t-\tau) K_m \big\|_2 \big\| v^{(q-1)} v_x (\tau) \big\|_1\; d\tau \\
&\leq  C\delta^q \int_0^t (1+t-\tau)^{-\frac{1}{2m}}(1+\tau)^{-\frac{q}{m}} \;d\tau 
\leq \delta^q (1+t)^{-\frac{1}{2m}}.
\end{align*}
In the ultima line we use that $q>m$ and \ref{gronw1}.\\
{\it{Proof of}} (\ref{svX2}),  by (\ref{sxk}), (\ref{gronw1}) and (\ref{svX0}) we have
\begin{align*}
\int_0^t \Big\| S_x(t-\tau) K_m \ast \big( v^q \big)_x (\tau) \Big\|_2\; d\tau 
\leq& \int_0^t \big\| S_x(t-\tau) K_m \big\|_2 \big\| v^{(q-1)} v_x (\tau) \big\|_1\; d\tau \\
\leq& \delta^q \int_0^t (1+t-\tau)^{-\frac{1}{2m}-\frac{1}{m}}(1+\tau)^{-\frac{q}{m}} \;d\tau 
\leq C\delta^q (1+t)^{-\frac{1}{2m}-\frac{1}{m}}.
\end{align*}
{\it{Proof of}} (\ref{svX3}),  by (\ref{ski}), (\ref{svX0}) and (\ref{gronw1}) we have
\begin{align*}
\int_0^t \Big\| S(t-\tau) K_m \ast \big( v^q \big)_x (\tau) \Big\|_{\infty}\; d\tau 
\leq& \int_0^t 
\big\| S(t-\tau) K_m \big\|_{\infty} \big\| v^{(q-1)} v_x (\tau) \big\|_1\; d\tau \\
&\leq \delta^q \int_0^t (1+t-\tau)^{-\frac{1}{m}}(1+\tau)^{-\frac{q}{m}} \;d\tau 
\leq C\delta^q (1+t)^{-\frac{1}{m}}.\quad\;{\BO}
\end{align*}
For the contraction property of $P$, we need the following Lema:
\begin{lem}\label{uvq}
Let $u(x,t), v(x,t)\in X_{\delta}$. Then there exist a positive constant $C$, independent
of $t$, such that
$$
\big\| \big( u^q \big)_x - \big( v^q\big)_x \big\|_1 
\leq C\delta^q M(u-v)(1+t)^{-q/m}, \qquad q>m.
$$
\end{lem}
\proof
We have that\;
$|u|^{q-1} u_x - |v|^{q-1} v_x = 
|u|^{q-1} \big( u_x - v_x \big) + v_x \big( |u|^{q-1}-  |v|^{q-1}\big)$, hence
\begin{align}\label{uvq1}
\big\| |u|^{q-1} u_x& - |v|^{q-1} v_x \big\|_1 \leq 
\big\| u^{q-1}\big\|_2 \;\big\| u_x - v_x  \big\|_2 
  + \big\| v_x \big\|_2 \;\big\| |u|^{q-1}-  |v|^{q-1} \big\|_2 \\
\leq& \big\| u\big\|_{\infty}^{q-2}\; \big\| u \big\|_2 \;\big\| u_x - v_x  \big\|_2
+ \big\| v_x \big\|_2 \;\bigg[ 
\big\| u-v \big\|_2 \Big( \big\|u \big\|_{\infty}^{q-1} + \big\|v \big\|_{\infty}^{q-1}\Big) 
\bigg] \;\;(\text{by}\;\;(\ref{holder1}))\nonumber \\
&\leq \delta^{q-1} \big\| u_x - v_x  \big\|_2\;(1+t)^{-\frac{q-2}{m} - \frac{1}{2m}}
+ \delta^q \big\| u-v \big\|_2 \;(1+t)^{-\frac{1}{2m}-\frac{1}{m}-\frac{q-1}{m}}\nonumber\\
&\quad =  \delta^{q-1} (1+t)^{\frac{1}{2m}+\frac{1}{m}}\; (1+t)^{-\frac{q}{m}} 
+\delta^q (1+t)^{\frac{1}{2m}}\big\| u-v \big\|_2 \;(1+t)^{-\frac{q}{m}-\frac{1}{m}}\nonumber\\
&\quad\quad \leq \delta^{q-1} M(u-v)(1+t)^{-\frac{q}{m}} + 
\delta^q M(u-v) (1+t)^{-\frac{q}{m}-\frac{1}{m}}\nonumber \\
&\quad\quad\quad \leq \delta^q M(u-v)(1+t)^{-\frac{q}{m}}.\hspace{6cm} {\BO}\nonumber
\end{align}
To prove Theorem 1, we need to prove that there exists the positive constant
$\delta$, such that the operator $P$ is a contraction mapping from 
$X_{\delta_1}$ into $X_{\delta_1}$.\\
{\it Step} 1.\;\;$P:X_{\delta} \rightarrow X_{\delta}$. For any 
$v_1(x,t) \in X_{\delta}$, and denoting $v=Pv_1$, we will prove that 
$v=Pv_1 \in X_{\delta}$ for some small $\delta >0$.\\
Indeed, using (\ref{sv}) and (\ref{svX1}) we have
\begin{align}\label{con1}
\big\| v(t) \big\|_2 =& \big\|P v_1(t) \big\|_2 \leq 
\big\| S(t)v_0 \big\|_2 + \int_0^t \big\| S(t-\tau)K_m \ast (v_1^q)_x (\tau) \big\|_2 d\tau  \\
&\leq C \big\| v_0 \big\|_{W^{1,1}} (1+t)^{-\frac{1}{2m}} +
C\delta^q (1+t)^{-\frac{1}{2m}}.\nonumber
\end{align}
Similarly, we have due to (\ref{sxv}) and (\ref{svX2})
\begin{gather}\label{con2}
\big\| v_x(t) \big\|_2 \leq 
C \big\| v_0 \big\|_{W^{2,1}}  (1+t)^{-\frac{1}{2m}-\frac{1}{m}} 
+ C\delta^q (1+t)^{-\frac{1}{2m}-\frac{1}{m} }. 
\end{gather}
By the same way, we can prove that
\begin{gather}\label{con3}
\big\| v(t) \big\|_{\infty} \leq 
C\big\| v_0 \big\|_{W^{2,1}}  (1+t)^{-\frac{1}{m}} 
+ C\delta^q (1+t)^{-\frac{1}{m} }. 
\end{gather}
Thus, combining (\ref{con1}), (\ref{con2}) and (\ref{con3}) implies that
$$
M(v)\leq C_1 \big( \big\| v_0 \big\|_{W^{2,1}} + \delta^q \big), \quad
q>m>3/2.
$$
Then there exist some small $\delta_2 >0$, such that $\delta_2^{q-1} < \frac{1}{2C_1}$.
Let  $ \big\| v_0 \big\|_{W^{2,1}} \leq \frac{\delta_2}{2C_1}$,
and $\delta \leq \delta_2$, then
$$
M(v)\leq C_1 \Big( \frac{\delta_2}{2C_1} + \delta_2^{q-1}\delta_2   \Big) <
\frac{\delta_2}{2} + \frac{\delta_2}{2} = \delta_2.
$$ 
We have proved $M(v)\leq \delta$ for some small $\delta $, namely, 
$P:X_{\delta} \rightarrow X_{\delta}$ for some small $\delta < \delta_2$.\\
{\it Step} 2.\;\; $P$ is a contraction in $X_{\delta}$.
Let $u,v \in X_{\delta}$, from Lema \ref{uvq} and Young inequality it follows that
\begin{align}\label{contra1}
\big\|Pu - Pv \big\|_2 &\leq \int_0^t \big\| S(t-\tau)K_m \big\|_2 \; 
\big\|(u^q)_x (\tau) - (v^q)_x (\tau) \big\|_1 d\tau  \\
&\leq C \delta^q M(u-v) \int_0^t (1+t-\tau)^{-\frac{1}{2m}}(1+\tau)^{-\frac{q}{m}} d\tau 
\leq C\delta^q M(u-v) (1+t)^{-\frac{1}{2m}}.\nonumber
\end{align}
We have, in the same way in (\ref{contra1})
\begin{multline}\label{contra2}
\big\|(Pu - Pv)_x \big\|_2 \leq \int_0^t \big\| S_x(t-\tau)K_m \big\|_2 \; 
\big\|(u^q)_x (\tau) - (v^q)_x (\tau) \big\|_1 d\tau  \\
\leq C \delta^q M(u-v) 
\int_0^t (1+t-\tau)^{-\frac{1}{2m}-\frac{1}{m}}\;(1+\tau)^{-\frac{q}{m}} d\tau 
\leq C\delta^q M(u-v) (1+t)^{-\frac{1}{2m} - \frac{1}{m}}.
\end{multline}
And also
\begin{align}\label{contra3}
\big\|Pu - Pv \big\|_{\infty} &\leq \int_0^t \big\| S(t-\tau)K_m \big\|_{\infty} \; 
\big\|(u^q)_x (\tau) - (v^q)_x (\tau) \big\|_1 d\tau  \\
&\leq C \delta^q M(u-v) \int_0^t (1+t-\tau)^{-\frac{1}{m}}(1+\tau)^{-\frac{q}{m}} d\tau 
\leq C\delta^q M(u-v) (1+t)^{-\frac{1}{m}}.\nonumber
\end{align}
Therefore, from (\ref{contra1}), (\ref{contra2}) and (\ref{contra3}), we obtain
\begin{gather}\label{contra4}
M(Pu-Pv)\leq C \delta^q M(u-v)
\end{gather}
Let choose $\delta \leq \delta_3 < \frac{1}{C^{1/q}}$; we have proved
$$
M(Pu-Pv)< M(u-v),
$$
i.e. $P:X_{\delta} \rightarrow X_{\delta}$ is a contraction for some small $\delta <\delta_3$.
Thank to steps 1 and 2, let $\delta_1 <\min \big\{ \delta_2, \delta_3 \big\}$, we have proved 
that the operator $P$ is contraction from $X_{\delta_1}$ to $X_{\delta_1}$. By the Banach's 
fixed point Theorem, we see that $P$ has a unique fixed point $v(x,t)$ in $X_{\delta_1}$.
This means that the integral equation (\ref{intequa}) has a unique global solution 
$v(x,t) \in X_{\delta_1}$. Thus, we have completed the proof of Theorem.
\hfill $\BO$ 
%
\setcounter{equation}{0}
\section{\hspace{-5mm}.\hspace{5mm} Asymptotic expansion: non-linear case}
In this section we prove some preliminary results  
which leads to the prove of the Theorem \ref{theo3}, i.e. to the 
asymptotic expansion of the solutions
of the equation (\ref{gene1}), where $q >m$, and $m>2$.\\
The solution of (\ref{gene1}) 
satisfies the integral equation (\ref{intequa})
obtained from the variation of constants formula. 
It is also convenient to recall that the solution of (\ref{gene1}) satisfy the 
decay properties (\ref{optimos}).\\ 
We now prove some preliminary results. \\
Noting that, from Theorem \ref{calorg}, with $N=0$ and 
${\cal{M}}=   \int_{\R}  u_0(x)dx$, we have
\begin{gather*}
t^{\frac{1}{m}(1-\frac{1}{p})}\big\| G_m(t)\ast u_0 - {\cal{M}}G_m(t) \big\|_p
\leq C t^{-\frac{1}{m}} \big\| |x| u_0 \big\|_1 \;,
\end{gather*}
for all $u_0 \in \Le^1((1+|x|)dx,\R)$.\\
Now, of this inequality, in view of the density of $\Le^1((1+|x|)dx,\R)$ in $\Le^1(\R)$, 
and in similar way to \cite{EsZu, Karch4} it is proven that:\; if $v_0 \in \Le^1(\R)$ then
\begin{gather}\label{prica}
t^{\frac{1}{m}(1-\frac{1}{p})} \big\| G_m(t)\ast u_0 - {\cal{M}}G_m(t) \big\|_p
\longrightarrow 0 \,,\quad\text{when}\quad t \rightarrow \infty
\end{gather}
Now, the decay rates from Lemmas \ref{lemaa} and \ref{lemab} are extend
to the case $p\in[2,\infty]$ , with $N=0$. Indeed, we have
%
%
\begin{lem}\label{casop}
There exist a constant $C=C(m)>0$, such that
\begin{align}
t^{\frac{1}{m}(1-\frac{1}{p})} 
\big\| S_{\varphi}(t) \ast v_0 (x)- G_m(t) \ast v_0 \big\|_p
\longrightarrow 0 \,,\quad\text{when}\quad t \rightarrow \infty, \label{casop1}\\
t^{\frac{1}{m}(1-\frac{1}{p})} 
\big\| S(t)v_0 (x)- S_{\varphi}(t) \ast v_0 (x) \big\|_p
\longrightarrow 0 \,,\quad\text{when}\quad t \rightarrow \infty, \label{casop2}
\end{align}
for all $p\in[2,\infty]$, and $v_0 \,\in \Le^1(\R) \cap \Le^2 (\R)$.
\end{lem}
%
\proof
For the proof we apply the interpolation inequality
$$
\| w\|_p \leq C \| w \|_2^{\frac{2}{p}} \| w\|_{\infty}^{1-\frac{2}{p}}\,,
\quad\forall\; p\in(2,\infty)\;.
$$
When $N=0$ in the Lemma \ref{lemaa}, we obtained the case $p=2$
\begin{gather*}
\nor S_{\varphi}(t)\ast u_0 -  G_m(t) \ast v_0 \nor_2 \leq 
C t^{ -1 - \frac{1}{2m} } \big\| u_0 \big\|_1. 
\end{gather*}
Now, we estimate the $\Le^{\infty}$-norm. As in (\ref{trans}) and (\ref{wplan2})
we decompose $G_m(x,t)= \int_{|\xi| \leq 1}...d\xi + \int_{|\xi| > 1}...d\xi$, then
\begin{align*}
\Big| S_{\varphi}(t)\ast u_0 -  G_m(t) \ast v_0 \Big| \leq C
\int_{|\xi| \leq 1} \Big| 
e^{-t|\xi|^m} \Big( e^{\frac{t|\xi|^{2m}}{1+|\xi|^m}} - 1 \Big) \Big| \;|\fo{v_0}|d\xi +
\int_{|\xi| > 1} e^{-t|\xi|^m} |\fo{v_0}|d\xi\,.
\end{align*}
Hence, of the Taylor expansion of the exponential function, $e^x -1\leq x e^x$, we have
\begin{align*}
\Big| S_{\varphi}(t)\ast u_0 -  G_m(t) \ast v_0 \Big| \leq C
\int_{|\xi| \leq 1} \Big| 
e^{-t\frac{|\xi|^m}{1+|xi|^m}} t|\xi|^{2m}\;|\fo{v_0}| d\xi +\int_{|\xi| > 1} e^{-t|\xi|^m} |\fo{v_0}|d\xi\,.
\end{align*}
If $|\xi|\leq 1$ then $\frac{|\xi|^m}{1+|\xi|^m} \geq \frac{|\xi|^m}{2}$ and if 
$|\xi| > 1$ then $2|\xi|^m > 1+|\xi|^m$, hence
\begin{align*}
\Big| S_{\varphi}(t)\ast u_0 -  G_m(t) \ast v_0 \Big| &\leq C
t\,\int_{|\xi| \leq 1}  
e^{-t\frac{|\xi|^m}{2}} |\xi|^{2m} d\xi \|v_0\|_1 + 
e^{-t} \int_{|\xi| > 1} e^{-t\frac{|\xi|^m}{2}} d\xi \|v_0\|_1\,.\\
&\leq C\Big(t^{-1 -\frac{1}{m}} + e^{-t}t^{-\frac{1}{m}}\Big)\|v_0\|_1 
\leq C t^{-1 -\frac{1}{m}} \|v_0\|_1.
\end{align*}
Hence, (\ref{casop1}) is a consequence of this inequality an interpolation formula.\\
The proff of (\ref{casop2}), is similar, using Lemma \ref{lemab}.\hfill$\BO$
%
In consequence, the first term of its asymptotic 
expansion of the linear solution $S(t)v_0$ is ${\cal{M}}G_m(t)$,
of course: 
\begin{theo}\label{prica1}
Let $v_0 \in \Le^1(\R) \cap \Le^2(\R)$. Then the linear solution $v(x,t)=S(t)v_0(x)$
of (\ref{genelin}) satisfies
$$
t^{\frac{1}{m}(1-\frac{1}{p})} \big\| S(t)v_0 (x)- {\cal{M}}G_m(t) \big\|_p
\longrightarrow 0 \,,\quad\text{when}\quad t \rightarrow \infty
$$
for all $p\in[2,\infty]$
\end{theo}
\proof
The relations (\ref{casop1}) and (\ref{casop2}) implies that
$$
t^{\frac{1}{m}(1-\frac{1}{p})} 
\big\| S(t)v_0 (x)-  G_m(t)\ast v_0 (x) \big\|_p
\longrightarrow 0 \,,\quad\text{when}\quad t \rightarrow \infty.
$$
This inequality and (\ref{prica}) provide the 
proof of Lemma.\hfill $\BO$
%
%
\begin{lem}\label{mite}
Se $a \in (-1,0]$ e $b\leq 0$, existe uma constante $C>0$,
independente de $t$, tal que
\begin{equation*}
\int_{0}^{t} (1+t-\tau)^{a}(1+\tau)^{b}d\tau \leq  \left\{
\begin{array}{l}
C(1+t)^{a +b +1}\quad \quad \quad for\quad b > -1 ; \\
C(1+t)^{a}\quad\quad\quad\quad \quad 
for\quad b < -1\; .
\end{array}
\right.
\end{equation*}
\end{lem}   
\begin{lem}\label{mita}
Assume $a \in (-1,0]$. There exist a constant $C$ independent of $t$ such that 
$$\int_0^t (1+t-\tau)^a (1+\tau)^{-1} d\tau \leq C(1+t)^{a} (1+log(1+t)).$$
\end{lem} 
%
\proof
After splitting the integral into $\int_0^{1/2} + \int_{1/2}^t$, the above inequality
is obtained by estimating each term by the supremum of one of the integrated factors.

The following Theorem say that if $ {\cal{M}}=\int u_0 (x) dx \neq 0$, then  
the first term of the asymptotic expansion of the solution $v(x,t)$ of 
(\ref{gene1}) is described by the 
fundamental solution of the linear equation (\ref{calge}).
\begin{theo}\label{primer}
Let $v=v(x,t)$ be the solution of (\ref{gene1}) corresponding to the initial 
data $v_0 \in \Le^1 (\R)\cap H^2 (\R)$. Then
\begin{align}
t^{\frac{1}{m}(1-\frac{1}{p})} 
\big\| v(t) -{\cal{M}}G_m(\cdot, t) \big\|_p & \leq \eta(t) \label{prim}\\
t^{\frac{1}{m}(q-\frac{1}{p})} 
\big\| v^q(t) -({\cal{M}}G_m(\cdot, t))^q \big\|_p & \leq \eta(t) \label{primA}.
\end{align}
such that $\lim\limits_{t\rightarrow \infty} \eta(t) =0$, 
with ${\cal{M}}=\int v_0(x) dx$ and $q>m$ 
\end{theo}
%
\proof
>From (\ref{intequa})
\begin{gather}\label{teo1}
\big\| v(t) - S(t)v_0 \big\|_p  
\leq \int_0^t \big\| S_x (t-\tau) K_m \big\|_p 
\big\| v(\tau) \big)\big\|_{\infty}^{q-m}  \big\| v(\tau) \big\|_m^m d\tau .
\end{gather}
>From (\ref{optimos}), by interpolation 
$\| v \|_m \leq C(1+t)^{-\frac{1}{m}(1-\frac{1}{m})}$. And by   
Lemma \ref{casp} following, we have
\begin{gather*}
\big\| v(t) - S(t)v_0 \big\|_p \leq 
C \int_0^t \big( 1+t-\tau \big)^{-\frac{1}{m}(1-\frac{1}{p}) - \frac{1}{m}} 
\big( 1+\tau\big)^{-\frac{1}{m}(q-1)} d\tau . 
\end{gather*}
Now, from Lemmas \ref{mita} and \ref{mite}, we have
\begin{gather*}
\big\| v(t) - S(t)v_0 \big\|_p \leq C \left\{
\begin{array}{lll}
t^{-\frac{1}{m}(1-\frac{1}{p}) + (\frac{m-q}{m})},\,\,m<q<m+1\\
t^{-\frac{1}{m}(1-\frac{1}{p}) -\frac{1}{m}} \log t,\,\, q=m+1\\
t^{-\frac{1}{m}(1-\frac{1}{p}) -\frac{1}{m}},\,q>m+1
\end{array}
\right.
\end{gather*}
Now, the proof of Theorem \ref{primer} is an immediate consequence of Theorema 
\ref{prica1}.\\
The proof of (\ref{primA}) result directly from (\ref{prim}) and (\ref{optimos}).
Indeed, it suffices to apply a simple consequence of the Lemma \ref{holder}
$$
\big\| v^q(t) -({\cal{M}}G_m(\cdot, t))^q \big\|_p 
\leq \| v(t) - {\cal{M}}G_m(\cdot, t)\|_p 
\Big( \big\|v(t) \big\|_{\infty}^{q-1} + \big\| {\cal{M}}G_m(\cdot, t) \big\|_{\infty}^{q-1}\Big)
\qquad\hfill{\BO}
$$
In the following Lemma, we extend the decay rates from Lemma \ref{skm} to the case
$p \in [0,\infty]$. For this, we use the interpolation inequality  
\begin{gather}\label{interpola}
\| g \|_{p} \leq C\| g \|_{\infty}^{1-\frac{1}{p}} \| g\|_1^{\frac{1}{p}}
\end{gather}
When $p=1$, observe that for all smooth rapidly decreasing functions $w=w(x)$ defined in $\R$,
\begin{gather}\label{norl1}
\| \fo{w} \|_{1} \leq C\| w \|_2^{1/2} \| \partial_x w \|_2^{1/2}.
\end{gather}
The proof of (\ref{norl1}) can be found e.g. in \cite{BDH} (example 2).
%
\begin{lem}\label{casp}
Let $m>2$.
Then, there exist a positive constant $C$, independent of $t$, such that
\begin{align}
\big\| S(t)K_m \big\|_p &
\leq C (1+t)^{-\frac{1}{m}(1-\frac{1}{p})} ,\label{casp1}\\
\big\| S_x(t)K_m \big\|_p &
\leq C (1+t)^{-\frac{1}{m}(1-\frac{1}{p}) -\frac{1}{m}}\,, \label{casp2}
\end{align}
for all $p \in [1, \infty]$.
\end{lem}
%
\proof
For the proof of (\ref{casp1}) we estimate the $\Le^1$-norm of $S(t)K_m$ 
and then we apply (\ref{interpola}), for this we use 
the inequality (\ref{norl1}) for the function 
$w(\xi) = \fo{K_m}e^{-t\Phi(\xi)}$, indeed we have
$$
\big| \partial_{\xi} w(\xi) \big|\leq C_m \Bigg[
t \Big(\frac{|\xi|^{m-1} + |\xi|^{m} + |\xi|^{2m}}{\big( 1+|\xi|^m \big)^3} \Big) +
\frac{|\xi|^{m-1}}{\big( 1+|\xi|^m \big)^2} \Bigg] 
e^{-\frac{t|\xi|^m}{1+|\xi|^m}}.
$$
Hence
\begin{align*}
\big| \partial_{\xi} w(\xi) \big| &
\leq C_m \big[ (t+1)|\xi|^{m-1} \big] e^{-\frac{t|\xi|^m}{2}}, 
\quad\text{if}\quad |\xi| \leq 1\,,\\
\big| \partial_{\xi} w(\xi) \big| & \leq C_m \Bigg[
\frac{t|\xi|^{2m}}{ \big( 1+|\xi|^m \big)^3} + 
\frac{|\xi|^{m-1}}{ \big( 1+|\xi|^m \big)^2} \Bigg] e^{-\frac{t}{2}}\,,
\quad\text{if}\quad |\xi| > 1\,.
\end{align*}
Then
\begin{multline*}
\|\partial_{\xi} w\|_2^2 =\int_{|\xi|\leq 1} + \int_{|\xi|>1}
\leq C (1+t)^2 \int_{|\xi|\leq 1} |\xi|^{2(m-1)} \big] e^{-t|\xi|^m} d\xi \\
+ C t^2 e^{-t}  \int_{|\xi|>1} \frac{|\xi|^{4m}}{ \big( 1+|\xi|^m \big)^6}d\xi
+C e^{-t}\int_{|\xi|>1} \frac{|\xi|^{2(m-1)}}{ \big( 1+|\xi|^m \big)^4} d\xi\,.
\end{multline*}
That is
\begin{gather}\label{wx2}
\|\partial_{\xi} w\|_2^2 \leq C (1+t)^{\frac{1}{m}} + Ct^2 e^{-t} + C e^{-t} \leq
C (1+t)^{\frac{1}{m}}.
\end{gather}
Now, substituting (\ref{sk}) and (\ref{wx2}) in (\ref{norl1}), we deduce that
\begin{gather}\label{jaja}
\|S(t)K_m\|_1 = \| \fo{w}\|_1 \leq C \Big( (1+t)^{-\frac{1}{2m}} \Big)^{1/2}
\Big( (1+t)^{\frac{1}{2m}} \Big)^{1/2} \leq C.
\end{gather}
Substituting (\ref{jaja}) and (\ref{ski}) in (\ref{interpola}) we have
\begin{gather*}
\| S(t)K_m \|_{p} \leq C(1+t)^{-\frac{1}{m}(1-\frac{1}{p})} \,.
\end{gather*}
The proof of (\ref{casp2}) is similar to (\ref{casp1}), in this case noting 
that, if $m>2$ 
\begin{gather*}
\big| S_x(t)K_m \big| \leq 
C\int \frac{ |\xi| \mbox{e}^{-\frac{t|\xi|^m}{1+|\xi|^m}} }{1+|\xi|^m} d\xi 
\leq C\int_0^1 |\xi| \mbox{e}^{-\frac{t|\xi|^m}{2}} d\xi +
C\int_1^{\infty} \frac{ |\xi| \mbox{e}^{-\frac{t}{2}} }{1+|\xi|^m}d\xi\\
\leq C(1+t)^{-\frac{1}{2m}} + Ce^{-\frac{t}{2}},
\end{gather*}
then
$$\big\| S_x(t)K_m \big\|_{\infty} \leq C(1+t)^{-\frac{1}{2m}}. \qquad{\BO}$$
\begin{remark}
When $m=2$, we have that $|\xi|/ (1+|\xi|^2) \not\in \Le^1(\R)$, however 
(\ref{casp2}) is valid for $p=\infty,\;\; m=2$, (see \cite{Karch1} Lemma 4.2) for
more details.
\end{remark}
%
\begin{lem}\label{caspg}
Let $G_m(x,t)$ the fundamental solution of the generalized heat equation 
(\ref{calge}). Then, for all $p\in[0,\infty]$,
\beqa\label{sg}
\big\| S(t)K_m - G_m(t) \big\|_p \leq Ct^{-\frac{1}{m}(1-\frac{1}{p})-\frac{1}{m}}\,,
\eeqa
\beqa\label{sderg}
\big\| \partial_x \big(S(t)K_m - G_m(t) \big) \big\|_p 
\leq Ct^{-\frac{1}{m}(1-\frac{1}{p})- \frac{2}{m}}\,.
\eeqa
\end{lem}
%
%
\proof
For the proof we use the same argument to that of Lemma \ref{casp}, 
we omit the details. \hfill ${\BO}$
\begin{remark}
When $p=2$, nothing that (\ref{sg}) is a consequence from the Theorem \ref{theo3},
indeed, if $v_0(x) = K_m(x)$, and as $\int K_m(x)dx =1$, 
whit $N=0$, we obtain the firs term in the asymptotic expansion i.e.
$$
 \big\| S(t)K_m - G_m(t) \big\|_2 \leq t^{-\frac{1}{2m}-\frac{1}{m}}\,,
$$
\end{remark}
Then as consequence of Lemmas \ref{casp} and \ref{caspg}  we have the following
Corollary:
%
\begin{cor}\label{kv}
Let $v(x,t)$ be the solution to (\ref{gene1}). Then there exists a constat $C>0$ 
such that for all $p\in [0,\infty]$ and $q \geq m+1$,
\begin{equation}\label{gkarch1}
\big\| S_x(t-\tau) K_m \ast v^q (\tau) \big\|_p \leq C \left\{
\begin{array}{l}
(t-\tau)^{-\frac{1}{m}(1-\frac{1}{p})-\frac{1}{m}} 
\tau^{- \frac{1}{m}(q-1)}\,; \\
(t-\tau)^{-\frac{1}{m}(1-\frac{1}{p})}
\tau^{-\frac{q}{m}}\,, 
\end{array}
\right.
\end{equation}
and
\begin{equation}\label{gkarch2}
\big\| \partial_x \big( S(t-\tau) K_m - G_m(t-\tau) \big) \ast v^q(\tau) \big\|_p 
\leq C \left\{
\begin{array}{l} 
(t-\tau)^{-\frac{1}{m}(1-\frac{1}{p})- \frac{2}{m}} \tau^{-\frac{1}{m}(q-1)}\,;\\
(t-\tau)^{-\frac{1}{m}(1-\frac{1}{p})-\frac{1}{m}} \tau^{- \frac{q}{m}}\,,
\end{array}
\right.
\end{equation}
for all $t>0$, $\tau \in(0,t)$.
\end{cor}
%
%
\proof
By interpolation it follows that
$$
\big\| v(t) \big\|_p \leq C (1+t)^{-\frac{1}{m}(1-\frac{1}{p})},
\quad\forall\;p \in [2,\infty]\,.
$$
And rememberig that 
$\big\| v_x(t) \big\|_2 \leq C (1+t)^{-\frac{1}{2m} - \frac{1}{m}}$.\\
These inequalities combined with Lemmas (\ref{casp}) and (\ref{caspg}) provide the 
proof of Corollary \ref{kv}. \hfill $\BO$\\
%
{\bf Proof of the Theorem \ref{theo3}}. 
We skip the proof of this Theorem because this differs from those Karch 
\cite{Karch4} in a few technical details, only.
%

%

\begin{thebibliography}{LLL}
%
\bibitem{ABFS}
L. Abdelouhab, J. L. Bona, M. Felland \& J. C. Saut ;
{\it Non local models for nonlinear dispersive waves}.
Physica, D40, 360-392 (1989).
%
\bibitem{ABS}
Ch. J. Amick, J. L. Bona \& M. E. Schonbek;
{\it Decay of solutions of some non linear wave equations}.
J. Differential Equations, 81, 1-49 (1989).
%
\bibitem{Benja1}
T. B. Benjamin;
{\it Lectures on Nonlinear Wave Motions}.
in Nonlinear Wave Motions, Alan Newell, ed., 
Lectures in Applied Mathematics No. 15, AMS. Providence, 
3-47 (1974).
%
\bibitem{Benja2}
T. B. Benjamin;
{\it A new kind of solitary wave}.
J. Fluid Mech. 245, 401-411 (1992).
%
\bibitem{Benja3}
T. B. Benjamin;
{\it Solitary and periodic waves of a new kind}.
Philos. Trans. Royal Soc. London A 354, 1775-1806 (1996).
%
\bibitem{BBM}
T. B. Benjamin, J. L. Bona \& J. J. Mahony;
{\it Model equations for long waves in nonlinear dispersive systems}.
Phil. Trans. Roy. Soc. London Ser. A, 272, 47-78 (1972).
%
\bibitem{Biler}
P. Biler;
{\it Asymptotic behavior in time of solutions to some equations
generalizing Korteweg-de Vries-Burger equations}.
Bull. Polish Acad. Sci. Ser. Math., 32, 275-282 (1984).
%
\bibitem{BDH}
P. Biler, J. Dziubanski \& W. Hebisch
{\it Scattering of small solutions to generalized Benjamin-Bona-Mahony equations
in several space dimensions}.
Comm. PDE, 17, 1737-1758 (1992)
%
\bibitem{BP}
V. Bisognin \& G. Perla Menzala;
{\it Decay rates of the solutions of nonlinear dispersive equations}.
Proc. Roy. Soc. Edinburgh Sect. A, 124, 1231-1246 (1994).
%
\bibitem{Bona1}
J. L. Bona;
{\it Solitary waves and other phenomena associated with model equations
for long waves}.
Fluid Dynamics Transaction, 10, 77-111 (1980).
%
\bibitem{Bona2}
J. L. Bona;
{\it On solitary waves and their role in the evolution of long waves}.
in Applications of Nonlinear Analysis in the Physical Sciences, H. Amann,
N. Bazley, K. Kirchg\"assner, ed. Pitman Press, London, 183-205 (1981)
%
\bibitem{BSht}
J. L. Bona \& R. Smith;
{\it The initial-value problem for the Korteweg-de Vries equation}.
Phil. Trans. Roy. Soc. London, Ser. A 278, 555-601 (1975).
%
\bibitem{Ca1}
A. Carpio;
{\it Asymptotic behavior for the vorticity equations in dimensions two
and three.}
Comm. Partial Differential Equations, 19, 827-872 (1994).
%
\bibitem{Ca2}
A. Carpio;
{\it Large-time behavior in incompressible Navier-Stokes equations}.
{SIAM J. Math. Anal., 27, 449-475 (1996).}
%
\bibitem{Cazweiss}
T. Cazenave \& F. B. Weissler;
{\it Asymptotically self-similar global solutions of the nonlinear
  Schr\"odinger and heat equations}.
{Math. Z., 228, No 1, 83-120 (1998).}
%
\bibitem{CBona}
H. Chen \& J. L. Bona;
{\it Existence and asymptotic proprties of solitary-wave solutions
of Benjamin-Type equations}.
{Preprint.}
%
\bibitem{Duozuazua}
J. Duoandikoetxea \& E. Zuazua;
{\it Moments, masses de Dirac et d\'ecomposition de fonctions}.
{C. R. Acad. Sci. Paris, t. 315, Serie I, 693-698 (1992).}
%
\bibitem{EsZu}
M. Escobedo \& E. Zuazua;
{\it Large time behavior for convection-diffusion equations in $\Rn$. }
J. Funct. Anal., 100, 119-161 (1991).
%
\bibitem{Friedman}
A. Friedman;
{\it Partial Diferential Equations}.
Rinehart and Winston, Inc. New York, 1969.
%
\bibitem{Henry}
D. B. Henry;
{\it How to remember the Sobolev inequalities}.
Lecture Notes in Mathematics, 957, 97-109, Springer, Berlin, 1982.
%
\bibitem{Karch1}
G. Karch;
{\it $L^p$-decay of solutions to dissipative-dispersive perturbations of
  conservation laws}.
Ann. Polon. Math., 67, 65-86 (1997).
%
\bibitem{Karch2}
G. Karch;
{\it Asymptotic behavior of solutions to some pseudo parabolic
  equations}.
Math. Methods Appl. Sci., 20, 271-289 (1997).
%
\bibitem{Karch4}
G. Karch;
{\it Large-time behavior of solutions to nonlinear wave equations:
  higher-order asymptotics}.
 { Math. Methods Appl. Sci.}, 22, 1671-1697, (1999).
%
\bibitem{Matsumura}
A. Matsumura;
{\it On the asymptotic behavior of solutions of semi-linear wave equation}.
 {Publ. RIMS. Kyoto Univ.}, 12, 169-189, (1976).
%
\bibitem{Mei}
M. Mei;
{\it Asymptotic expansion for Benjamin-Bona-Mahony-Burger equation
of generalized type}.
 {Nonlinear Analysis}, 33, 699-714, (1998).
%
\bibitem{Park}
M. A. Park;
{\it Large-time behavior of solution for 
Roseneau-Burger equations}.
 {Nonlinear Analysis}, 9, (1995).
%
\bibitem{razu}
R. Prado and E. Zuazua;
{\it Large-time behavior of solution for generalized 
Benjamin-Bona-Mahony-Burger equations}.
 {Reports LNCC}, 9, (2001).
%
\bibitem{Schonbek1}
M. E. Schonbek;
{\it Decay of solution to parabolic conservation laws}.
Comm. PDE., 7, 449-473 (1980).
%
\bibitem{Schonbek2}
M. E. Schonbek;
 {\it Uniform decay rates for parabolic conservation laws}.
 Nonlinear Analysis TMA, 10, 943-953 (1986).
%
\bibitem{Schonbek3}
M. E. Schonbek;
 {\it Large time behavior of solutions to the Navier-Stokes equations}.
Comm. PDE., 11, 733-763 (1986).
%
\bibitem{Segal}
I. E. Segal;
 {\it Quantization and dispersion for nonlinear relativistic equations}.
Mathematical Theory of Elementary Particles, MIT Press, Cambridge, MA. 79-108
(1966)
%
\bibitem{Stein}
E. M. Stein;
 {\it Singular Integrals and Differentiability Properties of Functions}.
{ Princeton University Press, Princeton: NJ (1970)}.
%
\bibitem{Ting}
T. W. Ting;
 {\it Certain non-steady flows of second order fluids}.
Arch. Rat. Mech. An., 18, 3-50 (1963).
%
\bibitem{Zhang2}
L. Zhang;
 {\it Decay of solutions of generalized Benjamin-Bona-Mahony-Burger
  equations in n-space dimensions}.
Nonlinear Analysis T.M.A., 25, 1343-1396 (1995).
%
\bibitem{Zuazua1}
E. Zuazua;
 {\it Weakly nonlinear large time behavior in scalar convection-diffusion
  equations}.
Differential Integral Equations, 6, 1481-1491 (1993).
%
\bibitem{Zuazua2}
E. Zuazua;
{\it Comportamiento asint\'otico de ecuaciones escalares de
  convecci\'on-difusi\'on}.
 Estudos e Comunica\c{c}\~oes do IMUFRJ 47, Rio de Janeiro, Brasil,
 1997.
\end{thebibliography}
\end{document}